\documentclass[11pt,reqno]{amsart}

\usepackage[T1]{fontenc}
\usepackage{lmodern}
\usepackage{amsmath,amssymb,amsthm,mathtools}
\usepackage{enumitem}
\usepackage{booktabs}
\usepackage{microtype}
\usepackage[hidelinks]{hyperref}
\usepackage[nameinlink,capitalize]{cleveref}
\setlist{nosep}
\numberwithin{equation}{section}

\newtheorem{theorem}{Theorem}[section]
\newtheorem{proposition}[theorem]{Proposition}
\newtheorem{lemma}[theorem]{Lemma}
\newtheorem{corollary}[theorem]{Corollary}

\newtheorem{branchinglemma}[theorem]{Lemma}
\crefname{branchinglemma}{Lemma}{Lemmas}
\Crefname{branchinglemma}{Lemma}{Lemmas}
\theoremstyle{definition}

\theoremstyle{remark}
\newtheorem{remark}[theorem]{Remark}

\DeclareMathOperator{\Ind}{Ind}
\DeclareMathOperator{\Hom}{Hom}
\DeclareMathOperator{\Res}{Res}

\DeclareMathOperator{\Sym}{Sym}
\DeclareMathOperator{\gr}{gr}

\newcommand{\C}{\mathbb C}
\newcommand{\R}{\mathbb R}
\newcommand{\g}{\mathfrak g}
\newcommand{\kfr}{\mathfrak k}
\newcommand{\pfr}{\mathfrak p}
\newcommand{\f}{\mathfrak f_4}
\newcommand{\s}{\mathfrak{sl}_2}
\newcommand{\PiMin}{\Pi_{\min}}
\newcommand{\ThetaB}{\Theta}
\newcommand{\thetaS}{\theta}
\newcommand{\boxt}{\boxtimes}

\title[The split exceptional theta correspondence]{The split exceptional theta correspondence for
\texorpdfstring{$\mathrm{PGL}_2\times F_4$}{PGL2 x F4}}
\author{Gordan Savin}

\begin{document}

\begin{abstract}
Let $G$ be the adjoint split real group of type $E_7$, and let $\PiMin$ be its
spherical minimal representation.  We study the exceptional theta correspondence
for the dual pair $\mathrm{PGL}_2(\R)\times F_{4(4)}$.
We prove that every irreducible $(\mathfrak{sl}_2,O(2))$-module has a
unique irreducible theta quotient, with one
exception: the small theta lift of $S^4(\C^2)$ is the sum of two irreducible
$F_{4(4)}$-modules obtained by restricting the minimal representation of
the split $E_6$. The proof is based on lifts of types of maximal compact subgroups.
These lifts are first computed on the graded module of the minimal representation,
the coordinate ring of the affine cone over the Grassmannian $\mathrm{Gr}(4,8)$, and then lifted
via the Poincar\'e--Birkhoff--Witt theorem to $\PiMin$. 
\end{abstract}

\maketitle
\tableofcontents

\section{Introduction}

This paper studies Howe duality, including its exceptional failure, for the
reductive dual pair
\[
             H\times F=\mathrm{PGL}_2(\R)\times F_{4(4)}
                  \subset E_{7(7)}^{\mathrm{ad}}.
\]
For the split ambient group, ``minimal representation'' means the
\emph{spherical} minimal representation.  Passing to the adjoint group is important.  Its identity
component has maximal compact group $\mathrm{SU}(8)/\mu_4$, but the reflection
which completes the circle $SO(2)$ to $O(2)$, the maximal compact subgroup of $\mathrm{PGL}_2(\R)$,
 is supplied by the disconnected normalizer.
It is this reflection that distinguishes the trivial and sign extensions of
the zero $SO(2)$-weight.

Let $\sigma$ be an irreducible admissible $(\s,O(2))$-module. Following \cite{BLS}
we define the big theta lift of $\sigma$ by
\[
\Theta_F(\sigma) =(\PiMin\otimes \sigma^{\vee})_{(\s,O(2))}.
\]
By \cite{LS}, $\Theta_F(\sigma)$ is a finite-length $(\f,K_F)$-module.
The main result is that $\Theta_F(\sigma)$ has a unique irreducible quotient,
with one exceptional failure.
The proof follows the ping--pong method of \cite{SavinPingPong},
developed there for the dual pair $\mathrm{SL}_2(\R)\times F_{4,1}$.
The main input to the ping--pong argument consists of the lifts of compact types.
More precisely, let $\tau$ be an $O(2)$-type. Define the lift of $\tau$ by
\[
\Theta_F(\tau)=(\PiMin\otimes \tau^{\vee})_{O(2)} \cong \Hom_{O(2)}( \tau, \Pi_{\min}).
\]
It is naturally a module for $U(\mathfrak f_4)$. We have a similar definition for the lift of
 $K_F=\mathrm{Sp}(3)\times\mathrm{SU}(2)$-types.

The main part of the paper is devoted to computing the lifts of compact types.
Let $\mathfrak g= \mathfrak k \oplus \mathfrak p$
be the Cartan decomposition for $E_{7(7)}^{\mathrm{ad}}$.
With $U=\C^8$, we have $\mathfrak p_{\mathbb C}=\Lambda^4U$.
The minimal representation has $K$-types that form a ladder, so that the corresponding
graded $S(\mathfrak p)$-module is
\[
R:=  \gr\PiMin=\bigoplus_{n\geq0}V_{\mathrm{SU}(8)}(n\varpi_4)
             =\C[\widehat{\mathrm{Gr}}(4,8)],
\]
the coordinate ring of the Pl\"ucker cone in $\Lambda^4 U$ \cite{BrylinskiKostant}.
All symmetric and enveloping algebras are over $\C$; in these algebras
and in ideals of $R$, complexification of real Lie algebras and Cartan
subspaces is understood.
This, via the filtered Nakayama lemma, turns the problem of computing lifts of
$O(2)$-types into computing the graded ring of coinvariants $R/\mathfrak p_F R$.
More precisely, this quotient is an $O(2)\times K_F$-module, with decomposition
\[
 R/\mathfrak p_F R
   =\bigoplus_{\tau\in\widehat{O(2)}}\tau\boxt\mathcal F_\tau,
 \qquad
 \mathcal F_\tau=\Hom_{O(2)}(\tau,R/\mathfrak p_F R).
\]
Our calculation determines these finite-dimensional graded $K_F$-modules.
Choosing homogeneous $K_F$-equivariant lifts of $\mathcal F_\tau$ to
$\Hom_{O(2)}(\tau,R)$ gives generators of this graded lift over
$S(\mathfrak p_F)$, by the graded Nakayama lemma.
Lifting once more through the filtration gives a copy of
$\mathcal F_\tau$ generating $\Theta_F(\tau)$ over $U(\mathfrak f_4)$.
In particular, $\Theta_F(\tau)$ is a quotient of the universal module
induced from $\mathcal F_\tau$:
\[
 U(\mathfrak f_4)\otimes_{U(\mathfrak k_F)}\mathcal F_\tau
                  \twoheadrightarrow\Theta_F(\tau).
\]

Our main result is most conveniently stated in terms of minimal $O(2)$-types.
Write
\[
 \tau_m=\C_m\oplus\C_{-m}\quad(m>0),\qquad
 \tau_0^+=\mathbf1,\quad \tau_0^-=\det
\]
for the irreducible representations of $O(2)$, with even indices in the
$\mathrm{PGL}_2$ normalization.

\begin{theorem}[Exceptional Howe correspondence]\label{thm:main}
Let $\sigma$ be an irreducible admissible $(\s,O(2))$-module belonging to
$\mathrm{PGL}_2(\R)$, and let $\ThetaB_F(\sigma)$ be its big theta lift in
$\PiMin$.
\begin{enumerate}
\item If $\sigma$ contains $\tau_0^+$ and
$\sigma\not\simeq S^4(\C^2)$, then $\ThetaB_F(\sigma)$ has exactly one
irreducible quotient, and that quotient is spherical.  For
$\sigma=S^4(\C^2)$ it has exactly two irreducible quotients; in the notation
of \cref{sec:exceptional-fibre},
\[
             \thetaS_F\bigl(S^4(\C^2)\bigr)
                  =\Pi_6^+\oplus\Pi_6^-.
\]
Here $\Pi_6^+$ is spherical and $\Pi_6^-$ has minimal $K_F$-type
$V_{\mathrm{Sp}(3)}(\varpi_2)\boxt\mathbf1$.
\item If $\sigma$ contains $\tau_0^-$, then $\ThetaB_F(\sigma)$ has exactly
one irreducible quotient, and that quotient contains
$V_{\mathrm{Sp}(3)}(\varpi_1)\boxt\C^2$ as its distinguished type.
\item If the minimal $O(2)$-type of $\sigma$ is $\tau_{2k}$, with $k\geq1$,
then
$\ThetaB_F(\sigma)$ has a unique irreducible quotient, containing
$V_{\mathrm{Sp}(3)}(k\varpi_2)\boxt\mathbf1$ with multiplicity one.
\item For the irreducible $(\s,O(2))$-module obtained by adjoining the
reflection to the lowest- and highest-weight modules with minimal type
$\tau_{2k}$, with $k\geq1$, this quotient is the module
$A_{\mathfrak q}((k-6)e_1)$ of \cref{sec:lowest-weight}.
\end{enumerate}
\end{theorem}

Thus a spherical representation of $\mathrm{PGL}_2(\R)$ lifts to a unique
irreducible spherical representation of $F_{4(4)}$, except at the single
infinitesimal character represented by $S^4(\C^2)$.  At that point there is
one additional nonspherical quotient.  The
global constant-term argument that produces the two quotients, and the Casimir
localization proving that there are no further exceptions, are given in
\cref{sec:exceptional-fibre}.

Section~\ref{sec:groups} constructs the dual pairs from the decomposition
$U=V\oplus W$, where $V$ and $W$ have dimensions $6$ and $2$.
We first introduce the dual pair $SO(2)\times E_{6(2)}$. Up to a harmless central extension,
 the maximal compact subgroup of the quasi-split $E_{6(2)}$ is
$K_E=\mathrm{SU}(V)\times\mathrm{SU}(W)$
 and the (complexified) space of non-compact directions is
 \[
 (\mathfrak p_E)_{\mathbb C}= \Lambda^3(V)\otimes W \subset \Lambda^4(U) = \mathfrak p_{\mathbb C}.
 \]
The group $E_{6(2)}$ has an involution $\iota$ whose fixed subgroup is $F_{4(4)}$.
In turn, the centralizer of $F_{4(4)}$ in the ambient group is
$H=\mathrm{PGL}_2(\R)$, whose maximal compact subgroup is covered by $O(2)$.

Section~\ref{sec:connected-lift} recalls Li's calculation of the lift from $SO(2)$ to $E_{6(2)}$.
A  calculation of the linear section $R/(\mathfrak p_E)$ of the Grassmannian cone gives a
direct proof that each nonzero lift is generated by its bottom
compact type. This, in turn, implies that the lift of any circle weight $2k$ is irreducible, including the small range
$-5\leq k\leq5$, completing Li's result.

Section~\ref{sec:forward-compact-lifts} computes the generators of the
$F_{4(4)}$-modules obtained by lifting $O(2)$-types, via a calculation of the linear section $R/(\mathfrak p_F)$.
Its nilradical is $(\mathfrak p_E)/(\mathfrak p_F)$ and it accounts for the additional
generators at the trivial and sign types.

Section~\ref{sec:reverse-compact-lifts} uses that $R$ is a unique factorization domain to prove
that it is free over $R^{K_F}=S(\mathfrak p_H)$.
The branching calculations then give $K_F$-isotypic spaces in $R$ as explicit, free $R^{K_F}$-modules.
This, via the PBW theorem, identifies the lifts of $K_F$-types in $\PiMin$ as
 modules induced from explicit finite-dimensional $O(2)$-representations.

Section~\ref{sec:pingpong} combines the two compact-type calculations
through see--saw reciprocity.  The resulting ping--pong argument proves
uniqueness of the irreducible theta quotient, with the spherical case
using the exceptional-fibre analysis of \cref{sec:exceptional-fibre}.

Frahm and Labriet \cite{FrahmLabriet} determine the unitary decomposition
of minimal representations under a family of dual pairs with a
$\mathrm{PGL}_2$ member, including
$\mathrm{PGL}_2(\R)\times F_{4(4)}$ in split $E_7$.
Their Theorems~A and~B(2) identify the discrete summands as tensor
products of discrete series of $\mathrm{PGL}_2(\R)$ and
cohomologically induced representations of $F_{4(4)}$.
Section~\ref{sec:lowest-weight} matches their parameters with ours and
restricts the discrete projections to Harish--Chandra modules.
Together with the uniqueness in part~(3), this proves part~(4).

Finally, Section~\ref{sec:E6-F4-restriction} proves that the spherical
minimal representation of split $E_6$ restricts to two irreducible
$F_{4(4)}$-modules.  The quotient of its Lagrangian Grassmannian cone ring
by $\mathfrak p_F$ is the square-zero algebra
$\C\oplus\Lambda^2_0\C^6$.  Its two generating types, separated by the
diagram involution, give the two constituents of the exceptional fibre.

A $p$-adic version of Theorem \ref{thm:main} was established in \cite{KarasiewiczSavin}.
It is interesting to note that it also features a similar exceptional fibre.
For a global variant of this correspondence, see \cite{Shan}.

\paragraph{\bf Acknowledgments.}
I gratefully acknowledge support from the Simons Foundation through gift No.~946504,
and thank the University of Utah for providing access to OpenAI. 
ChatGTP-6 Astra was used for assistance with typesetting, 
and to compute linear section of Grassmannians. 
The proofs included here are self contained. 
The linear section $R/(\mathfrak p_E)$ is trivial to describe, directly from Pl\"ucker relations. 
On the the other hand, due to the presence of nilpotent elements, $R/(\mathfrak p_F)$ 
is more difficult to compute. In the end, its nilpotent radical is nicely described
 in terms of an apolar Artin algebra corresponding to $V\otimes W$ \cite{Weyman}. 

\section{The groups, Cartan decompositions, and the see--saw}
\label{sec:groups}

\subsection{The split E7 model}

Let $\g=\mathfrak e_{7(7)}$ and choose a Cartan decomposition
$\g=\kfr\oplus\pfr$.  With $U=\C^8$, its complexification is
\begin{equation}\label{eq:e7-cartan}
 \kfr_\C=\mathfrak{sl}(U),\qquad
 \pfr_\C=\Lambda^4U,\qquad \kfr=\mathfrak{su}(U).
\end{equation}
The bracket on $\pfr_\C$ is the equivariant map
$\Lambda^4U\otimes\Lambda^4U\to\mathfrak{sl}(U)$ determined up to scale
by a volume form on $U$.
We use $\mathrm{SU}(U)$ as a compact covering group for the calculations.
Group products below are understood up to finite central intersections;
the ambient real group is the split adjoint group together with the
normalizer elements used below.

\subsection{The circle and E6}

Choose an orthogonal decomposition
\begin{equation}\label{eq:block}
 U=V\oplus W,\qquad \dim_\C V=6,\quad\dim_\C W=2,
\end{equation}
and put $K_E=\mathrm{SU}(V)\times \mathrm{SU}(W)\subset \mathrm{SU}(U)$.
Its connected centralizer is the block-scalar circle.
In the normalization appropriate to $\mathrm{PGL}_2$, the infinitesimal
generator of the circle group on $U$ is
\[
 S=\frac{i}{2}\operatorname{diag}(I_6,-3I_2).
\]
The kernel of the action of $\mathrm{SU}(U)$ on
$\mathfrak{sl}(U)\oplus\Lambda^4U$ is $\mu_4$: a central scalar
$\zeta I_U$ acts on the second summand by $\zeta^4$.
Thus $\exp(\pi S)=iI_U$ is already trivial in the adjoint group.
Throughout, $SO(2)$ denotes the standard compact double cover of this
circle, with parameter $t$ modulo $2\pi$; its kernel is $\{0,\pi\}$.
Likewise $O(2)$ denotes the corresponding cover of the compact group
$PO(2)$ in $\mathrm{PGL}_2(\R)$.  All representations under consideration
are trivial on the kernel, which explains the even character weights.  The exterior-power decomposition is
\begin{equation}\label{eq:p-decomp-first}
 \pfr_\C=\Lambda^4V
 \oplus(\Lambda^3V\otimes W)
 \oplus(\Lambda^2V\otimes\Lambda^2W),
 \qquad\text{weights }2,\ 0,\ -2.
\end{equation}
Thus the zero-weight space is
\[
 (\pfr_E)_\C=\Lambda^3V\otimes W,
 \qquad
 \kfr_E=\mathfrak{su}(V)\oplus\mathfrak{su}(W).
\]
Together these form the real Lie algebra $\mathfrak e_{6(2)}$ inside
$\g$.  This is Li's quasi-split real form $E_{6,4}$, of real rank four.
Indeed, the circle centralizer has Cartan decomposition with compact part
$\R S\oplus\kfr_E$ and noncompact part $\pfr_E$; its derived algebra is
$\mathfrak e_{6(2)}$.  In particular the correct centralizer statements are
\[
 \mathfrak z_{\g}(\R S)=\R S\oplus\mathfrak e_{6(2)},\qquad
 \mathfrak z_{\g}(\mathfrak e_{6(2)})=\R S,
\]
or, at the connected group level,
\[
 C_G(SO(2))^\circ=SO(2)\cdot E_{6(2)},\qquad
 C_G(E_{6(2)})^\circ=SO(2).
\]
This constructs the commuting pair $SO(2)\times E_{6(2)}$ before
introducing its normalizer.

\subsection{The reflection}

We now define the reflection $r$ which extends the circle to $O(2)$.
Its conjugation action on $\mathrm{SU}(U)/\mu_4$ will be the involution
$\sigma$ below, and its action on $\pfr_\C$ will be the complex-linear
involution $\mathcal R$.

Choose an ordered orthonormal basis
\[
 (e_1,f_1,e_2,f_2,e_3,f_3,e_4,f_4)
\]
adapted to $V\oplus W$, and set
\[
 \mathrm{vol}=e_1\wedge f_1\wedge e_2\wedge f_2
                 \wedge e_3\wedge f_3\wedge e_4\wedge f_4.
\]
Put
\[
 w=\begin{pmatrix}0&-1\\1&0\end{pmatrix},\qquad
 J=\operatorname{diag}(w,w,w,w).
\]
The quaternionic structure on $U$ is the antilinear map
$\mathcal J(u)=J\overline u$.  It satisfies
$\mathcal J^2=-1$ and $\mathcal J(zu)=\overline z\,\mathcal J(u)$.
The induced involution on $\mathrm{SU}(U)$ is
\begin{equation}\label{eq:reflection-compact}
 \sigma(g)=J\overline gJ^{-1}.
\end{equation}
Its fixed subgroup is $\mathrm{Sp}(4)$, the unitary group preserving this
quaternionic structure.  It preserves $K_E$, with fixed subgroup
\[
 K_F=K_E^\sigma=\mathrm{Sp}(3)\times \mathrm{SU}(2),
\]
and it sends $S$ to $-S$.

To specify the action of $r$ on $\pfr_\C$, let $u_I$ be the wedge of
four distinct basis vectors in the chosen order, and let $u_{I^c}$ be
the complementary basis wedge, also in the chosen order.  Define
$*_0(u_I)$ to be $u_{I^c}$ with the sign chosen so that
\[
                         u_I\wedge *_0(u_I)=\mathrm{vol},
\]
and extend $*_0$ complex-linearly to $\Lambda^4U$.
Thus the sign is simply the sign needed to reorder the concatenated
lists $I,I^c$ into the chosen order of all eight basis vectors.
This is volume-form duality followed by the basis identification
$U^*\simeq U$; there is no complex conjugation in $*_0$.
Define
\begin{equation}\label{eq:reflection-exterior}
 \mathcal R=(\Lambda^4J)\circ *_0.
\end{equation}
Then
\[
 \mathcal R^2=1,\qquad
 \mathcal R\,\Lambda^4g=\Lambda^4(\sigma(g))\,\mathcal R
 \quad(g\in \mathrm{SU}(U)).
\]
The differential of $\sigma$ on $\mathfrak{sl}(U)$ is
$X\mapsto-JX^{\mathsf t}J^{-1}$.  Together with $\mathcal R$ it preserves
the brackets in \cref{eq:e7-cartan} and the split real form.
For example, the standard real structure on $\pfr_\C$ is
$u\mapsto *_0\overline u$, and it commutes with $\mathcal R$.
These actions define the reflection $r$ in the adjoint real group:
conjugation by $r$ on $\mathrm{SU}(U)/\mu_4$ is induced by $\sigma$,
and $\operatorname{Ad}(r)|_{\pfr_\C}=\mathcal R$.
It inverts $SO(2)$, so $\langle SO(2),r\rangle\simeq O(2)$ in our
normalization.

For the signs on the middle summand, put
\[
 \Omega_V=\sum_{i=1}^3e_i\wedge f_i,\qquad
 \Omega_W=e_4\wedge f_4.
\]
Symplectic contraction gives
\begin{equation}\label{eq:sympl-exterior}
 \Lambda^2V=\mathbf1\oplus V(\varpi_2),\qquad
 \Lambda^3V=\Lambda^3_0V\oplus(\Omega_V\wedge V)
           =V(\varpi_3)\oplus V(\varpi_1).
\end{equation}
Directly from \cref{eq:reflection-exterior},
\begin{align*}
 \mathcal R(e_1\wedge e_2\wedge e_3\wedge e_4)
   &=e_1\wedge e_2\wedge e_3\wedge e_4,\\
 \mathcal R(\Omega_V\wedge e_1\wedge e_4)
   &=-\Omega_V\wedge e_1\wedge e_4.
\end{align*}
Since $\mathcal R$ commutes with $K_F$, these calculations determine its
action on the two irreducible summands of $\Lambda^3V\otimes W$:
it is $+1$ on $\Lambda^3_0V\otimes W$ and $-1$ on
$(\Omega_V\wedge V)\otimes W$.
It also exchanges the weight-$2$ and weight-$-2$ summands of
\cref{eq:p-decomp-first}.  Consequently the decomposition under
$K_F\times O(2)$ is
\begin{equation}\label{eq:70}
\begin{split}
 \Lambda^4U={}&
 \bigl(V(\varpi_3)\boxt W\bigr)\boxt\tau_0^+\\
 &\oplus\bigl(\mathbf1\oplus V(\varpi_2)\bigr)\boxt\tau_2\\
 &\oplus\bigl(V(\varpi_1)\boxt W\bigr)\boxt\tau_0^-.
\end{split}
\end{equation}
Its dimensions are $28+30+12=70$.
This derives the sign of the last summand, which is invisible to $SO(2)$.

\subsection{F4 and the rank-one centralizer}

The $O(2)$-fixed part of \cref{eq:70} is
\[
 (\pfr_F)_\C=V_{\mathrm{Sp}(3)}(\varpi_3)\boxt W.
\]
In the compact algebra, the $O(2)$-centralizer is $\kfr_F$:
the circle centralizer is $\R S\oplus\kfr_E$, while $r$ negates $S$
and fixes $\kfr_F$ inside $\kfr_E$.
Thus $\kfr_F\oplus\pfr_F$ is the real Lie algebra
$\mathfrak f_{4(4)}$, the fixed algebra of $r$ inside $\mathfrak e_{6(2)}$.

Conversely, the $K_F$-fixed subspace of $\pfr_\C$ is
\begin{equation}\label{eq:pH-generators}
 (\pfr_H)_\C=\C z_+\oplus\C z_-,
 \qquad z_+=\tfrac12\Omega_V^2,\quad
 z_-=\Omega_V\wedge\Omega_W.
\end{equation}
The normalization is chosen because
$\mathcal R(\Omega_V^2)=2\Omega_V\wedge\Omega_W$;
hence $r$ interchanges $z_+$ and $z_-$.
These vectors commute also with $\pfr_F$.  Indeed their brackets with
$\pfr_F$ would lie in the weight-$\pm2$ subspaces of $\kfr_\C$,
each of $K_F$-type $V(\varpi_1)\boxt W$, whereas
$\pfr_F$ has type $V(\varpi_3)\boxt W$; equivariance forces those
brackets to vanish.  Thus $(\pfr_H)_\C$ is precisely the centralizer
of $\mathfrak f_{4(4)}$ in $\pfr_\C$.
Together with $\R S$, its real form gives the centralizing
$\mathfrak{sl}_2(\R)$, and adjoining $r$ gives the rank-one group
$\mathrm{PGL}_2(\R)$ in the adjoint normalization.

We therefore obtain the see--saw
\begin{equation}\label{eq:seesaw}
\begin{matrix}
 \mathrm{PGL}_2(\R)&\quad\times\quad&F_{4(4)}\\[2mm]
 \cup&&\cap\\[-1mm]
 SO(2)&\quad\times\quad&E_{6(2)}.
\end{matrix}
\end{equation}
The reflection normalizes, rather than centralizes, the $E_6$ factor:
\begin{equation}\label{eq:normalizer-action}
 rtr^{-1}=t^{-1},\qquad rgr^{-1}=\iota(g)
 \quad(t\in SO(2),\ g\in E_{6(2)}),
\end{equation}
where $\iota$ is the outer involution with fixed Lie algebra
$\mathfrak f_{4(4)}$.
If $\widetilde M$ is generated by $SO(2)E_{6(2)}$ and $r$, then, up to
finite central identifications,
\begin{equation}\label{eq:normalizer-quotients}
 \widetilde M/E_{6(2)}\simeq O(2),\qquad
 \widetilde M/SO(2)\simeq E_{6(2)}\rtimes\langle\iota\rangle.
\end{equation}
In particular the connected centralizer of $O(2)$ is $F_{4(4)}$.
The groups $O(2)$ and $F_{4(4)}$ commute, whereas $O(2)$ and $E_{6(2)}$
do not.

\section{The lift from \texorpdfstring{$SO(2)$}{SO(2)} to \texorpdfstring{$E_{6(2)}$}{E6(2)}}
\label{sec:connected-lift}

We recall Li's calculation of the lifts from $SO(2)$ to $E_{6(2)}$
and give a Pl\"ucker-relations proof that each nonzero lift is generated
by a single $K_E$-type.  Unitarity then proves irreducibility, including
the small range $0\leq k<6$ not covered by Li's cohomologically induced
identification.

\subsection{The ladder model}

Let $\PiMin$ be the spherical minimal Harish--Chandra module of
$E_{7(7)}^{\mathrm{ad}}$.  Its $K$-type decomposition is multiplicity free:
\begin{equation}\label{eq:ladder}
       \PiMin\big|_{\mathrm{SU}(8)}
           =\bigoplus_{n\geq0}V_{\mathrm{SU}(8)}(n\varpi_4).
\end{equation}
The natural good filtration has associated graded algebra
\begin{equation}\label{eq:grass-ring}
 R:=\gr\PiMin=\bigoplus_{n\geq0}V_{\mathrm{SU}(8)}(n\varpi_4)
       \simeq\C[\widehat{\mathrm{Gr}}(4,U)],
\end{equation}
as in \cite{BrylinskiKostant}.  This is the homogeneous coordinate ring of
the affine Pl\"ucker cone; see, for example, \cite{Weyman}.

For the character $\chi_m$ of $SO(2)$ of weight $m$, the connected theta lift
\begin{equation}\label{eq:connected-E6-lift}
 \Sigma_m=\ThetaB_{E_6}(\chi_m)
       :=\Hom_{SO(2)}(\chi_m,\PiMin)
\end{equation}
is an $E_{6(2)}$-module.  Only even weights $m$ occur.  We write
$m=2k$ or $m=-2k$, where $k=|m|/2\geq0$; the two choices coincide
when $m=k=0$.  This is the correspondence calculated by Li.

\subsection{Li's calculation}

Write $\omega_i$ for the fundamental weights of $\mathrm{SU}(6)$.
Duality interchanges $\omega_2$ and $\omega_4$.
In our normalization, Li's integer $k$ corresponds to circle weight $2k$.
He proves that
$\sigma_k^+=\Sigma_{2k}$ has
\begin{equation}\label{eq:li-ktypes}
 \sigma_k^+\big|_{\mathrm{SU}(6)\times\mathrm{SU}(2)}
 =\bigoplus_{n,r\geq0}
 V_{\mathrm{SU}(6)}\bigl(n\omega_2+r\omega_3+(n+k)\omega_4\bigr)
       \boxt\Sym^rW,
\end{equation}
and $\sigma_k^-=\Sigma_{-2k}$ is given by the contragredient formula
\cite{Li}.  For $k\geq6$ he identifies
\begin{equation}\label{eq:li-aq}
       \sigma_k^\pm\simeq A_{\mathfrak q^\pm}
                    \bigl((k-6)\omega^\pm\bigr).
\end{equation}
A generalized Blattner formula gives the $K$-types.  A Vogan--Zuckerman
lemma \cite{VoganZuckerman} then identifies the irreducible unitary module
from its lowest $K$-type and infinitesimal character.

\subsection{A Pl\"ucker proof of finite generation}

The decomposition in \cref{eq:p-decomp-first} will be written
\begin{equation}\label{eq:ABC-decomposition}
 \Lambda^4U=A\oplus B\oplus C,
 \qquad
 \begin{cases}
  A=\Lambda^4V,\\
  B=\Lambda^3V\otimes W,\\
  C=\Lambda^2V\otimes\Lambda^2W.
 \end{cases}
\end{equation}
The three summands have $SO(2)$-weights $2,0,-2$.  In particular,
$B=(\pfr_E)_\C$ is the complexified noncompact tangent space of
$E_{6(2)}$.  Let $J=(B)\subset R$.

\begin{proposition}[The reduced linear section]
\label[proposition]{prop:B-linear-section}
The ideal $J$ is radical, and
\begin{equation}\label{eq:reduced-ring}
 R/J\simeq R^+\mathop{\times}_{\C}R^-
      =\C\oplus R^+_{>0}\oplus R^-_{>0},
 \qquad
 \begin{cases}
  R^+=\C[\widehat{\mathrm{Gr}}(4,V)],\\
  R^-=\C[\widehat{\mathrm{Gr}}(2,V)].
 \end{cases}
\end{equation}
Consequently the projective zero locus of $B$ on
$\mathrm{Gr}(4,V\oplus W)$ is the disjoint union
\begin{equation}\label{eq:two-components}
 \mathrm{Gr}(4,V)\ \sqcup\
 \{L\oplus W:L\in\mathrm{Gr}(2,V)\}.
\end{equation}
\end{proposition}

\begin{proof}
We start with the standard Pl\"ucker relation, recorded below as
\eqref{eq:standard-pluecker-relation}:
\[
 \sum_{a=1}^5(-1)^{a-1}
 p_{u_1u_2u_3t_a}\,
 p_{t_1\cdots\widehat{t_a}\cdots t_5}=0.
\]
Here all $u_i,t_j$ are vectors in $U$, and $p$ is multilinear and
alternating in its four vector arguments.  If all these vectors belong
to $V$, we obtain the Pl\"ucker relations for $\mathrm{Gr}(4,V)$
among the $A$-coordinates.

Let $x,y$ be a basis of $W$, and identify $C$ with $\Lambda^2V$ using
$x\wedge y$.  Set $u_2=t_4=x$ and $u_3=t_5=y$, with
$u_1,t_1,t_2,t_3\in V$.  The last two terms vanish by alternation,
and the relation becomes
\[
 p_{u_1t_1xy}p_{t_2t_3xy}
 -p_{u_1t_2xy}p_{t_1t_3xy}
 +p_{u_1t_3xy}p_{t_1t_2xy}=0.
\]
These are the Pl\"ucker relations for $\mathrm{Gr}(2,V)$ among the
$C$-coordinates.

Finally, set $u_3=x$ and $t_5=y$, with
$u_1,u_2,t_1,t_2,t_3,t_4\in V$.  Moving $x$ to the last position in
each of the first four factors gives
\begin{align*}
 p_{t_1t_2t_3t_4}\,p_{u_1u_2xy}
  ={}&p_{u_1u_2t_1x}\,p_{t_2t_3t_4y}
      -p_{u_1u_2t_2x}\,p_{t_1t_3t_4y}\\
    &+p_{u_1u_2t_3x}\,p_{t_1t_2t_4y}
      -p_{u_1u_2t_4x}\,p_{t_1t_2t_3y}.
\end{align*}
Every factor on the right belongs to $B$.  As the vectors vary, the
two factors on the left span $A$ and $C$, so $AC\subset B^2$ in $R$.
Setting $B=0$ therefore gives $AC=0$ in $R/J$.

Modulo $B$, only the $A$- and $C$-coordinates remain.  Together with
their pure Pl\"ucker relations, the vanishing of $AC$ gives
\cref{eq:reduced-ring}: the two Grassmannian cones lie in the linear
section, so there are no further relations.  Both component rings are
domains, and the fibre product embeds in $R^+\times R^-$; hence $R/J$
is reduced.  Its two affine cones meet only at their common vertex,
and projectivizing gives \cref{eq:two-components}.
\end{proof}

Borel--Weil identifies the degree-$k$ parts of the two Grassmannian
cone rings in \cref{eq:reduced-ring} as
\[
 R^+_k=V_{\mathrm{SU}(6)}(k\omega_4)\boxt\mathbf1,
 \qquad
 R^-_k=V_{\mathrm{SU}(6)}(k\omega_2)\boxt\mathbf1,
\]
where the second factor is the trivial representation of
$\mathrm{SU}(W)$.  Since $A$ and $C$ have circle weights $2$ and $-2$,
these components have weights $2k$ and $-2k$, respectively.
Consequently, as an $SO(2)\times K_E$-module,
\begin{equation}\label{eq:BW-two-cones}
\begin{split}
 R/J={}&\chi_0\boxt\mathbf1\\
 &\oplus\bigoplus_{k\geq1}
 \chi_{2k}\boxt\bigl(V_{\mathrm{SU}(6)}(k\omega_4)\boxt\mathbf1\bigr)\\
 &\oplus\bigoplus_{k\geq1}
 \chi_{-2k}\boxt\bigl(V_{\mathrm{SU}(6)}(k\omega_2)\boxt\mathbf1\bigr).
\end{split}
\end{equation}
Each summand indexed by $k$ lies in degree $k$, and the two constants
are identified in degree zero.  Thus, for each even circle weight $m$,
the quotient has exactly one $K_E$-type, in degree $|m|/2$.

\begin{corollary}[The lift from $SO(2)$]\label[corollary]{cor:connected-lifts}
For every $k\geq0$, the $E_{6(2)}$-module $\Sigma_{2k}$ is generated by
its single $K_E$-type
$V_{\mathrm{SU}(6)}(k\omega_4)\boxt\mathbf1$, and $\Sigma_{-2k}$ is
generated by $V_{\mathrm{SU}(6)}(k\omega_2)\boxt\mathbf1$.  Every nonzero
connected theta lift $\Sigma_m$ occurring in $\PiMin$ is irreducible.
\end{corollary}

\begin{proof}
Equip $\Sigma_m$ with the filtration induced by the good filtration of
$\PiMin$ in \cref{eq:grass-ring}.  Taking $SO(2)$-isotypic components
commutes with passage to the associated graded module, so
\[
 \gr\Sigma_m=\Hom_{SO(2)}(\chi_m,R).
\]
The action of $B=(\pfr_E)_\C$ has symbol multiplication by $B$ in $R$.
Moreover, $B$ has circle weight zero, so exactness of taking isotypic
components gives
\[
 (\gr\Sigma_m)/B(\gr\Sigma_m)
       \simeq\Hom_{SO(2)}(\chi_m,R/J).
\]
By \cref{eq:BW-two-cones}, the right-hand side is the claimed generating
$K_E$-type in degree $k=|m|/2$: it is
$V_{\mathrm{SU}(6)}(k\omega_4)\boxt\mathbf1$ for $m=2k$,
$V_{\mathrm{SU}(6)}(k\omega_2)\boxt\mathbf1$ for $m=-2k$,
and the trivial type for $m=0$.
The graded Nakayama lemma now shows that $\gr\Sigma_m$ is generated
over $\Sym(B)$ by this type.  Choose a $K_E$-equivariant lift of it
to $\Sigma_m$.  By PBW, its $U(\mathfrak e_{6,\C})$-span has
associated graded module equal to $\gr\Sigma_m$, and hence equals
$\Sigma_m$.  This is the filtered Nakayama argument proving the
asserted generation.

Each $\Sigma_m$ is an admissible unitary representation, being an
$SO(2)$-isotypic summand of the unitary representation $\PiMin$.  Its
generating $K_E$-type is irreducible and occurs with multiplicity one.
If $\Sigma_m$ had two nonzero irreducible summands, orthogonal complete
reducibility would place the generating type in only one of them, contrary
to generation.  Thus $\Sigma_m$ is irreducible.  Combined with
\cref{eq:li-ktypes}, this completes the connected correspondence for every
$k$; for $k\geq6$ it agrees with Li's cohomologically induced
identification \cref{eq:li-aq}.
\end{proof}

\section{The \texorpdfstring{$O(2)$}{O(2)}-multiplicity modules for \texorpdfstring{$F_{4(4)}$}{F4(4)}}
\label{sec:forward-compact-lifts}

For an irreducible $O(2)$-type $\tau$, define the multiplicity module
\begin{equation}\label{eq:O2-F-module}
             \mathcal A_\tau=\Theta_F(\tau)=\Hom_{O(2)}(\tau,\PiMin).
\end{equation}
Because the centralizer of the full $O(2)$ is $F_{4(4)}$, this is an
$F_{4(4)}$-module.  We determine finite-dimensional $K_F$-modules
that generate these lifts.

\subsection{The normalizer refinement}

For $m>0$, restriction from $O(2)$ to $SO(2)$ and the action of $r$ give
\begin{equation}\label{eq:O2-from-SO2}
 \mathcal A_{\tau_m}
      \simeq \Sigma_m\big|_{F_{4(4)}}
      \simeq \Sigma_{-m}\big|_{F_{4(4)}}.
\end{equation}
More precisely, $r$ exchanges the two connected lifts $\Sigma_m$ and
$\Sigma_{-m}$ while implementing $\iota$ on $E_{6(2)}$; their direct sum is
the $\tau_m$-isotypic summand of $\PiMin$.  In particular, for the
normalization used by Li, $r$ exchanges $\sigma_k^+=\Sigma_{2k}$ and
$\sigma_k^-=\Sigma_{-2k}$.

At weight zero the reflection acts within $\Sigma_0$, and its two
eigenspaces are precisely the multiplicity modules for the two extensions
of the trivial $SO(2)$-character:
\begin{equation}\label{eq:zero-weight-splitting}
       \Sigma_0\big|_{F_{4(4)}}
          =\mathcal A_{\tau_0^+}\oplus\mathcal A_{\tau_0^-}.
\end{equation}
The signs in this decomposition are fixed by \cref{eq:70}: the
$V(\varpi_3)\boxt W$ direction is even and the
$V(\varpi_1)\boxt W$ direction is odd.

We now prove finite generation of all the modules $\mathcal A_\tau$.
For $\tau=\tau_{2k}$ with $k>0$, this refines the calculation of
\cref{cor:connected-lifts}; at weight zero the nilpotent structure of the
linear section distinguishes the trivial and sign eigenspaces.

\subsection{The bottom layer}

The symplectic branching rule for the rectangular representation is
particularly simple.

\begin{proposition}[Rectangular branching]\label[proposition]{prop:rectangular}
For every $k\geq0$,
\begin{equation}\label{eq:rectangular}
 \Res^{\mathrm{SU}(6)}_{\mathrm{Sp}(3)}V_{\mathrm{SU}(6)}(k\omega_4)
      =\bigoplus_{j=0}^k V_{\mathrm{Sp}(3)}(j\varpi_2),
\end{equation}
with multiplicity one.
\end{proposition}

\begin{proof}
Identify $V(k\omega_4)$ with the Schur module of rectangular partition
$(k,k,k,k)$.  Littlewood's restriction rule from $\mathrm{GL}(6,\C)$ to $\mathrm{Sp}(6,\C)$
\cite{KoikeTerada} says that
the multiplicity of the symplectic module $V_{\mathrm{Sp}(3)}(\mu)$ is
$\sum_\delta c^{(k^4)}_{\mu,(2\delta)' }$.  A rectangle admits such an even
column complement only when $\mu=(j,j)$, $0\leq j\leq k$, and the relevant
Littlewood--Richardson coefficient is one.  The partition $(j,j)$ is the
weight $j\varpi_2$.  No other $\mu$ occurs.
\end{proof}

For reference, the first cases are
\[
\begin{array}{c|c|c}
k&\Res V(k\omega_4)&\text{dimension check}\\ \midrule
1&\mathbf1\oplus V(\varpi_2)&15=1+14\\
2&\mathbf1\oplus V(\varpi_2)\oplus V(2\varpi_2)&105=1+14+90\\
3&\mathbf1\oplus V(\varpi_2)\oplus V(2\varpi_2)\oplus
V(3\varpi_2)&490=1+14+90+385.
\end{array}
\]

Let
\begin{equation}\label{eq:Fk}
             \mathcal F_k=\bigoplus_{j=0}^k
                 V_{\mathrm{Sp}(3)}(j\varpi_2)\boxt\mathbf1.
\end{equation}
It is the restriction of the bottom $K_E$-type in
\cref{eq:li-ktypes}.

\subsection{The ring obtained by killing the noncompact directions}

Put
\[
 P=V(\varpi_3)\boxt W,\qquad Q=V(\varpi_1)\boxt W,
 \qquad B=P\oplus Q,
\]
and let $I=(P)\subset R$.  Thus $P=(\mathfrak p_F)_\C$, and our goal
is to determine the ring $R/I$.  Its degree-one generators are $A,C,Q$.
Recall that $J=(B)=(P,Q)$ and that $R/J$ was determined in
\cref{prop:B-linear-section}.

Choose symplectic basis vectors $e_1,e_2,e_3,f_1,f_2,f_3$ of $V$, with
$\langle e_i,f_j\rangle=\delta_{ij}$, and a basis $x,y$ of $W$.
The quotient $B/P\simeq V\otimes W$ is given by contraction
\[
 \gamma(u\wedge v\wedge w)
 =\langle u,v\rangle w-\langle u,w\rangle v+\langle v,w\rangle u.
\]
Normalize the generators $q_{v,z}$ of $Q$ so that, modulo $I$,
\[
                   p_{uvwz}=q_{\gamma(u\wedge v\wedge w),z}
                   \qquad(u,v,w\in V,\ z\in W).
\]
The notation is linear in each vector; all Pl\"ucker coordinates are
alternating in their indices.  The twelve generators $q_{e_i,x}$,
$q_{f_i,x}$, $q_{e_i,y}$, $q_{f_i,y}$ form a $6\times2$ matrix.

\paragraph{The infinitesimal scheme in the $Q$-directions.}
Before determining the mixed products, consider the graded algebra
$\mathcal D$ generated by $Q=V\otimes W$ with relations
\begin{equation}\label{eq:ring-quadratic-relations}
 \begin{gathered}
 q_{u,x}q_{v,x}=q_{u,y}q_{v,y}=0,\qquad
 q_{u,x}q_{v,y}+q_{v,x}q_{u,y}=0\quad(u,v\in V),\\
                         \sum_{i=1}^3q_{e_i,x}q_{f_i,y}=0.
 \end{gathered}
\end{equation}
The mixed products are alternating in $u,v$, and the last relation kills
their symplectic trace.  Put $\Omega=\sum_i e_i\wedge f_i$ and
$T=\Lambda^2V/\C\Omega$.  Then
\begin{equation}\label{eq:Q-algebra-pieces}
 \mathcal D_0=\C,\qquad \mathcal D_1=Q,\qquad
 \mathcal D_2=T,\qquad \mathcal D_d=0\quad(d\geq3),
\end{equation}
where $q_{u,x}q_{v,y}$ corresponds to the class of $u\wedge v$.
Indeed, every cubic monomial has two factors from the same column, so
it is zero already by the first relations.  Thus
$\mathcal D_{>0}^3=0$, and the Hilbert series is $1+12t+14t^2$.

This is a scheme of length $27$ supported at the origin, not an
additional reduced cone: every matrix entry has square zero and therefore
vanishes at every complex point.  To compare with the Segre cone, recall
\[
 \Sym^2(V\otimes W)
 =\bigl(\Sym^2V\otimes\Sym^2W\bigr)
   \oplus\bigl(\Lambda^2V\otimes\Lambda^2W\bigr).
\]
The rank-one, or Segre, cone is obtained by killing the second summand,
the $2\times2$ minors.  Here the first summand is killed, whereas the
minors survive modulo their symplectic trace.  They are nilpotent
functions describing infinitesimal structure at the origin.  The trace
relation uses the symplectic form and is invariant under
$\mathrm{Sp}(3)\times\mathrm{SU}(2)$; before imposing it, the quadratic
relations are invariant under $\mathrm{SL}(V)\times\mathrm{SL}(W)$.

\begin{lemma}[The quadratic relations in $Q$]\label[lemma]{lem:Q-relations}
The generators of $Q$ in $R/I$ satisfy
\eqref{eq:ring-quadratic-relations}.  Consequently multiplication induces
an equivariant surjection
\[
 T=V(\varpi_2)\boxt\mathbf1\longrightarrow Q^2,
\]
and $Q^3=0$.  In particular, $Q^2$ is either zero or isomorphic to $T$.
\end{lemma}

\begin{proof}
We use the standard Pl\"ucker relations
\begin{equation}\label{eq:standard-pluecker-relation}
 \sum_{a=1}^5(-1)^{a-1}
 p_{s_1s_2s_3t_a}\,
 p_{t_1\cdots\widehat{t_a}\cdots t_5}=0.
\end{equation}

The choice
\[
 (s_1,s_2,s_3)=(e_1,e_2,x),\qquad
 (t_1,\ldots,t_5)=(e_1,e_3,f_2,f_3,x)
\]
has just one nonzero term after contraction, namely $-q_{e_1,x}^2$.
Changing symplectic bases and polarizing in $V$ and $W$ gives the first
two relations in \eqref{eq:ring-quadratic-relations}.

For the trace relation, specialize to
\[
 (s_1,s_2,s_3)=(e_1,e_2,e_3),\qquad
 (t_1,\ldots,t_5)=(f_1,f_2,f_3,x,y).
\]
The two $BB$-terms vanish because both triples are isotropic, leaving
\[
 p_{e_1e_2e_3f_1}p_{f_2f_3xy}
 -p_{e_1e_2e_3f_2}p_{f_1f_3xy}
 +p_{e_1e_2e_3f_3}p_{f_1f_2xy}=0.
\]
Put $t_i=q_{e_i,x}q_{f_i,y}$.  Applying
\eqref{eq:standard-pluecker-relation} with triple $(f_i,f_j,x)$ and
quintuple $(e_1,e_2,e_3,f_k,y)$, for
\[
                 (i,j,k)=(2,3,1),\ (1,3,2),\ (1,2,3),
\]
respectively, evaluates the three
products as $-t_2-t_3$, $t_1+t_3$, and $-t_1-t_2$.
Their alternating sum is $-2(t_1+t_2+t_3)$, proving the trace relation.
The description of $\mathcal D$ above now gives the asserted surjection
and $Q^3=0$.
\end{proof}

\begin{lemma}[The mixed products in $R/I$]\label[lemma]{lem:nilpotent-quotient}
In $R/I$, the subspaces spanned by products satisfy
\[
                   AQ=CQ=0,\qquad AC=Q^2\simeq
                   V(\varpi_2)\boxt\mathbf1.
\]
\end{lemma}

\begin{proof}
\emph{The products $CQ$ and $AQ$.}
Take the triple $(u,x,y)$ and the quintuple $(v_1,v_2,v_3,v_4,x)$
in \eqref{eq:standard-pluecker-relation}.  Up to an overall sign, the result is
\[
 \sum_{a=1}^4(-1)^{a-1}p_{uv_axy}\,
 q_{\gamma(v_1\wedge\cdots\widehat{v_a}\cdots\wedge v_4),x}=0.
\]
Take $(v_1,v_2,v_3,v_4)=(e_i,e_j,e_k,f_k)$, where $i,j,k$ are
distinct.  Only two terms remain:
\[
                 p_{ue_ixy}q_{e_j,x}-p_{ue_jxy}q_{e_i,x}=0.
\]
The same calculation with $(f_i,f_j,e_k,f_k)$ or
$(e_i,f_j,e_k,f_k)$ allows us to switch $f_i,f_j$ or $e_i,f_j$,
respectively, when $i\ne j$.

For the remaining pairs, take
$(v_1,v_2,v_3,v_4)=(e_i,f_i,e_j,f_j)$.  The relation becomes
\[
 \bigl(p_{ue_ixy}q_{f_i,x}-p_{uf_ixy}q_{e_i,x}\bigr)
 +\bigl(p_{ue_jxy}q_{f_j,x}-p_{uf_jxy}q_{e_j,x}\bigr)=0.
\]
This holds for each of the pairs $(1,2)$, $(1,3)$, and $(2,3)$.
Adding the first two relations and subtracting the third shows that
twice the first parenthesis for $i=1$ is zero; the other two cases
follow in the same way.  Thus we may also switch $e_i,f_i$.
By bilinearity, for arbitrary $v,w\in V$ we have
\[
                       p_{uvxy}q_{w,x}=p_{uwxy}q_{v,x}.
\]
Antisymmetry in $u,v$ now gives
\begin{align*}
 p_{uvxy}q_{w,x}
 &=p_{uwxy}q_{v,x}=-p_{wuxy}q_{v,x}\\
 &=-p_{wvxy}q_{u,x}=p_{vwxy}q_{u,x}=-p_{uvxy}q_{w,x},
\end{align*}
so this product is zero.  Interchanging $x,y$ gives $CQ=0$.
The reflection $r$ interchanges $A$ and $C$ and preserves $I$ and $Q$;
therefore $AQ=0$ as well.

\emph{The products $AC$.}
The specialization
\[
 (s_1,s_2,s_3)=(u,v,x),\qquad
 (t_1,\ldots,t_5)=(v_1,v_2,v_3,v_4,y)
\]
is the one used in \cref{prop:B-linear-section}, with each
$B$-coordinate $p_{abc z}$ now replaced by
$q_{\gamma(a\wedge b\wedge c),z}$.  Thus $AC\subset B^2$ in $R$
becomes $AC\subset Q^2$ in $R/I$.  One particular case is
\[
              p_{e_1e_2f_1f_2}p_{e_1e_2xy}
                     =2q_{e_1,x}q_{e_2,y}.
\]
By \cref{lem:Q-relations}, it remains to check that $Q^2\ne0$:
then this product is nonzero, and equivariance and the irreducibility
of $T$ give $AC=Q^2=T$.

\emph{Nonvanishing in degree two.}
The branching calculation in \cref{lem:covariant-branching}(1),
with $k=1$, gives a copy of $T=V(\varpi_2)\boxt\mathbf1$ in
$R_2^{O(2)}$.  By \cref{eq:70}, $P$ is fixed pointwise by $O(2)$
and $R_1^{O(2)}=P$.  Exactness of taking invariants therefore gives
\[
 (R/I)_2^{O(2)}\simeq
 R_2^{O(2)}\big/
 \operatorname{im}\bigl(\Sym^2P\longrightarrow R_2\bigr).
\]
On the other hand, since $P=V(\varpi_3)\boxt W$,
\[
 (\Sym^2P)^{\mathrm{SU}(2)}
       =\Lambda^2V(\varpi_3)
       =\mathbf1\oplus V(2\varpi_2).
\]
Thus $\Sym^2P$ contains no $T$, and the copy of $T$ in
$R_2^{O(2)}$ survives in $(R/I)_2^{O(2)}$.
The generators $A,C,Q$ have circle weights $2,-2,0$, respectively,
and $O(2)$ acts on $Q$ by the sign character.  Hence the degree-two
$O(2)$-invariants in $R/I$ are spanned by invariant $AC$-products
and $Q^2$.  Since $AC\subset Q^2$, we obtain
\[
                         (R/I)_2^{O(2)}=Q^2\ne0.
\]
The surjection $T\twoheadrightarrow Q^2$ of \cref{lem:Q-relations}
is consequently an isomorphism.  The mixed product displayed above
then gives $AC=Q^2\simeq T$, as required.
\end{proof}

The consequences needed for theta lifting are now immediate.  The ideal
$N=J/I$ is generated by $Q$, has $N_1=Q$, $N_2=Q^2$, and $N_d=0$
for $d\geq3$.  Its quotient is $R/J$, and the positive-degree parts of
both pure component rings annihilate $N$.  Thus $N$ is the nilradical
and $\sqrt I=J$.  In the representation notation used elsewhere in the
paper, this says
\begin{equation}\label{eq:nilpotent-pieces}
 \begin{gathered}
 N_1=\bigl(V(\varpi_1)\boxt W\bigr)\boxt\tau_0^-,\qquad
 N_2=\bigl(V(\varpi_2)\boxt\mathbf1\bigr)\boxt\tau_0^+,\\
                         N_d=0\quad(d\geq3).
 \end{gathered}
\end{equation}
The signs follow from \cref{eq:70}.

\begin{theorem}[Finite generation of the $O(2)$ multiplicities]
\label{thm:finite-generation}
For $k>0$, the $F_{4(4)}$-module $\mathcal A_{\tau_{2k}}$ is generated over
$U(\f)$ by the finite-dimensional $K_F$-module $\mathcal F_k$ of
\cref{eq:Fk}.  The zero-type multiplicity modules have generators
\begin{align*}
 \mathcal A_{\tau_0^+} &: \quad \mathbf1\oplus
               V(\varpi_2)\boxt\mathbf1,\\
 \mathcal A_{\tau_0^-} &: \quad V(\varpi_1)\boxt W.
\end{align*}
In particular all the modules $\mathcal A_\tau$ are finitely generated
$U(\f)$-modules.
\end{theorem}

\begin{proof}
Give $\mathcal A_\tau$ the filtration induced from \cref{eq:grass-ring}.
The action of $\mathfrak p_F$ has symbol multiplication by $P$.  Therefore
the quotient by $\mathfrak p_F$ is controlled by $R/I$.  The
$SO(2)\times K_E$-decomposition of its reduction $R/J$ is given in
\cref{eq:BW-two-cones}.  For $k>0$, the weight-$2k$ and weight-$-2k$
parts are $V_{\mathrm{SU}(6)}(k\omega_4)\boxt\mathbf1$ and its dual
$V_{\mathrm{SU}(6)}(k\omega_2)\boxt\mathbf1$, respectively.
Restricting either module to $K_F$ and applying \cref{prop:rectangular}
gives exactly $\mathcal F_k$; the reflection in $O(2)$ exchanges the
two cones.

Finally, the nilradical $J/I$ described in
\cref{lem:nilpotent-quotient} has $O(2)$-types only $\tau_0^-$
in degree one and $\tau_0^+$ in degree two.  It therefore contributes
nothing for $k>0$ and gives precisely the two displayed corrections at
weight zero.  Lifting a basis of the coinvariant quotient through the good
filtration and applying the filtered Nakayama lemma proves generation over
$U(\f)$.
\end{proof}

\section{Reverse compact-type lifts to \texorpdfstring{$(\mathfrak{sl}_2,O(2))$}{(sl2,O(2))}}
\label{sec:reverse-compact-lifts}

For a $K_F$-type $\kappa$, put
\[
 \mathcal M(\kappa)=\Hom_{K_F}(\kappa,\PiMin).
\]
The centralizer of $F$ acts on this space, making it a
$(\s,O(2))$-module.  We shall compute this module by considering the graded version of this problem.

\subsection{Freeness}
 Recall that $R=\bigoplus_{n\geq0}R_n$, where
$R_n=V_{\mathrm{SU}(8)}(n\varpi_4)$.  Use the forms
$\omega=\Omega_V$ and $\epsilon_W=\Omega_W$ fixed in Section~\ref{sec:groups}.  Thus
the generators of \cref{eq:pH-generators} are
\[
 z_+=\tfrac12\omega\wedge\omega\in\Lambda^4V,
 \qquad
 z_-=\omega\wedge\epsilon_W
       \in\Lambda^2V\otimes\Lambda^2W
\]
and are $K_F$-fixed vectors of $SO(2)$-weights $2$ and $-2$ in $R_1$.
The reflection $r$ already fixed in \cref{eq:reflection-exterior}
satisfies $r z_+=z_-$ and $r z_-=z_+$.  Put
$Z=\C z_+\oplus\C z_-$.
Put $X=\C_2\oplus\C_{-2}\simeq Z$, the complexified noncompact tangent
representation of $\mathrm{PGL}_2/O(2)$.
The associated graded of the enveloping-algebra action on the spherical
vector gives the invariant-algebra isomorphism \cite{LS}
\[
\Sym(X) \simeq R^{K_F}=\C[z_+,z_-].
\]

\begin{proposition}\label[proposition]{prop:A-freeness} The ring $R$ is a free graded $R^{K_F}$-module.
\end{proposition}
We use the following elementary graded criterion.
\begin{lemma}\label[lemma]{lemma:freeness}
Let $M$ be a graded $\C[z_+,z_-]$-module, bounded below and with
finite-dimensional graded pieces, where $\deg z_+=\deg z_-=1$.
Assume that multiplication by
$z_+$ on $M$ is injective, and multiplication by $z_-$ on $M/z_+M$ is also injective. Then $M$ is
a free graded $\C[z_+,z_-]$-module.
\end{lemma}
\begin{proof}
Put $\overline M=M/(z_+M + z_-M)$.  A graded
linear section $\overline M\to M$ induces, by multiplication and
induction on degree (using the lower bound), a surjective homorphism of graded $\C[z_+,z_-]$-modules
\[
                  N:= \C[z_+,z_-]\otimes_{\C}\overline M\rightarrow M.
\]
Let $H_M(t)$ denote the Poincar\'e series of the graded module $M$.
The two injectivity assumptions imply
\[
H_{M/z_+ M}(t)=(1-t)H_M(t) \quad\text{and}\quad H_{\overline M}(t)=(1-t) H_{M/z_+ M}(t).
\]
Hence
\[
H_M(t)= \frac{H_{\overline M}(t)}{(1-t)^2}= H_N(t),
\]
and the surjection is an isomorphism.
\end{proof}

\begin{proof}[Proof of \cref{prop:A-freeness}]
We show that $R$ satisfies the conditions of the lemma.
Firstly, $R$ is an integral domain, so multiplication by $z_+$ is injective.
Secondly, it is well known that the Pl\"ucker ring $R$ is
a unique factorization domain.  One way to see this is that the
Grassmannian has Picard group generated by the Pl\"ucker hyperplane class;
the normal affine cone therefore has trivial divisor class group \cite{Popov}.
Since $R$ is connected and positively graded, every nonzero element of
$R_1$ is irreducible.  In particular $z_+$ is prime and $R/(z_+)$ is a
domain.  This argument does not require $z_+$ to be a Pl\"ucker coordinate.
The class of $z_-$ in this quotient is nonzero, since the two degree-one
vectors have different $SO(2)$-weights and so $z_-$ is not in $(z_+)$.
It follows that $z_+,z_-$ satisfy the conditions of Lemma \ref{lemma:freeness}.
\end{proof}

For any $K_F$-type $\kappa$, write
\[
 M_\kappa=\Hom_{K_F}(\kappa,R).
\]
The $K_F$-isotypic projections commute with multiplication by $R^{K_F}$.
The regular sequence $z_+,z_-$ on $R$ therefore remains regular on
each $M_\kappa$, so Lemma~\ref{lemma:freeness} makes each multiplicity
space graded free as well.  We now describe these modules for the
$\kappa$ appearing in the lifts of $O(2)$-types.

\begin{branchinglemma}[The covariant branching algebra]
\label[lemma]{lem:covariant-branching}
The $R^{K_F}$-modules $M_\kappa$ have the following form.
\begin{enumerate}
\item For $T_k=V_{\mathrm{Sp}(3)}(k\varpi_2)\boxt\mathbf1$, multiplication
induces an isomorphism of graded $O(2)$-modules
\begin{align*}
 R^{K_F}\otimes\Hom_{K_F}(T_k,R_k)
       &\xrightarrow{\ \sim\ }
          M_{T_k}\\
 \Hom_{K_F}(T_k,R_k)&\simeq\Sym^k Z.
\end{align*}
Thus the $T_k$-covariants form a free module of rank $k+1$ over
$R^{K_F}$, generated in degree $k$.
\item The $Q_F=V(\varpi_1)\boxt W$ covariants form a free rank-one module
over $R^{K_F}$, generated in degree one by the sign line:
\[
 R^{K_F}\otimes \det[1]
      \simeq M_{Q_F}.
\]
Here $[1]$ records that the generator has degree one.
\end{enumerate}
\end{branchinglemma}

\begin{proof}
The freeness argument above shows that $M_\kappa$ is a free $R^{K_F}$-module.
We proceed to compute the Poincar\'e series of $M_\kappa$.
To that end we need to decompose each $R_n$,
which we view as a representation of $\mathrm{GL}(U)$
of highest weight $(n,n,n,n,0,0,0,0)$.  Then, as a
$\mathrm{GL}(V)\times\mathrm{GL}(W)$-module,
\begin{equation}\label{eq:Cauchy-rectangle}
 R_n =\bigoplus_{0\leq s\leq r\leq n}
       V(\lambda_{r,s})\boxt V((n-s,n-r))
\end{equation}
where $\lambda_{r,s}=(n,n,r,s,0,0)$.
The $\mathrm{SL}(W)$-fixed summands are precisely those with
$r=s$, $0\leq r\leq n$.  The corresponding $SO(2)$-weight is
\begin{equation}\label{eq:weight-nr}
 \tfrac12|\lambda_{r,r}|-\tfrac32|(n-r,n-r)|=4r-2n.
\end{equation}

Next, we restrict to $\mathrm{Sp}(V)$.  For part~(1) we need the
highest weight $(k,k,0)$, or $k\varpi_2$ in fundamental-weight
coordinates; for part~(2) we need the standard representation of highest
weight $(1,0,0)$.  Both multiplicities follow directly from the
Littlewood--Richardson restriction rule \cite{KoikeTerada}.

\begin{lemma}\label[lemma]{lem:paired-branching}
For $p,r,k\geq0$, the multiplicity of the $\mathrm{Sp}(3)$-module
of highest weight $(k,k,0)$ in the restriction of the
$\mathrm{SU}(6)$-module of highest weight $(p+r,p+r,r,r,0,0)$ is
\begin{align*}
 &\#\{h\in\mathbb Z:\max(0,k-p)\leq h\leq\min(k,r)\}\\
 &\hspace{3em}=\max\{0,\ 1+\min(k,p,r,p+r-k)\}.
\end{align*}
In particular it is nonzero exactly when $p+r\geq k$.

For $1\leq r\leq n$, the $\mathrm{Sp}(3)$-module of highest weight
$(1,0,0)$ occurs with multiplicity one in the restriction of the
$\mathrm{U}(6)$-module of highest weight $(n,n,r,r-1,0,0)$.
\end{lemma}

\begin{proof}
Put $N=p+r$.  If $N<k$, the diagram $(k,k)$ is not contained in
$(N,N,r,r)$ and the multiplicity is zero.  Otherwise put $L=N-k$.
The restriction rule counts Littlewood--Richardson tableaux of shape
$(N,N,r,r)/(k,k)$ whose content has even columns.  Such content has
the form $(x,x,y,y)$.  The lattice-word and column conditions force
the rows, from left to right, to be
\[
 1^L,\qquad 2^L,\qquad
 1^h3^{r-h},\qquad 2^h4^{r-h}.
\]
Column strictness gives $0\leq h\leq\min(k,r)$, and the lattice-word
condition gives $r-h\leq L$, or $h\geq k-p$.  Conversely these
inequalities make the displayed filling a valid tableau, of content
$(L+h,L+h,r-h,r-h)$.  Distinct $h$ give distinct contents, each with
multiplicity one.  No finite-rank correction contributes: the possible
length-four universal symplectic characters specialize to zero in rank
three by the symplectic specialization rule of \cite{KoikeTerada}.

For the second assertion put $\lambda=(n,n,r,r-1)$.
The multiplicity of $(1,0,0)$ is
\[
 \sum_{\eta\ \text{with even columns}}c_{(1),\eta}^{\lambda}
 =\sum_{\eta\ \text{with even columns}}c_{\eta,(1)}^{\lambda}.
\]
By the one-box Littlewood--Richardson rule, a coefficient in the last
sum is one precisely when $\lambda/\eta$ is a single box, and is zero
otherwise.  The only removal leaving even columns is the last box of
the third row, giving $\eta=(n,n,r-1,r-1)$.  The single box, filled
with $1$, is the unique tableau, proving multiplicity one.
\end{proof}

Apply the lemma with $p=n-r$ and sum over $0\leq r\leq n$ in
\cref{eq:Cauchy-rectangle}.  The admissible pairs are
$0\leq h\leq k$ and $h\leq r\leq n-k+h$.
Thus $M_{T_k}$ is zero in degrees $n<k$, and its degree-$n$ dimension
is $(k+1)(n-k+1)$ for $n\geq k$.  Hence
\[
 H_{M_{T_k}}(t)=
 \frac{(k+1)t^k}
 {(1-t)^2}.
\]
Thus $\overline M_{T_k}:=M_{T_k}/(z_+,z_-)M_{T_k}$ is concentrated
in degree $k$.
Its degree-$k$ part is $\Hom_{K_F}(T_k,R_k)$.  Moreover,
\cref{eq:weight-nr} gives
\[
 \left.\Hom_{K_F}(T_k,R_k)\right|_{SO(2)}
       \simeq\bigoplus_{j=0}^k\C_{2k-4j}.
\]
This determines the circle action, but for even $k$ it does not yet
determine the reflection on the zero-weight line.

We check the reflection using multiplication in $R$.  Fix a Borel
subgroup of $K_{F,\C}$ and identify the multiplicity of $T_k$ with the
space of highest vectors of weight $k\varpi_2$, trivial on the
$\mathrm{SU}(2)$ factor.  By \cref{eq:70}, there is a nonzero such
vector $u_+\in R_1$ of highest weight $\varpi_2$ and circle weight $2$.
Put $u_-=r u_+$.  Since $r$ centralizes $K_F$ and $r^2=1$, the vector
$u_-$ has the same highest weight, has circle weight $-2$, and satisfies
$r u_-=u_+$.  The Cartan-product algebra $R$ carries the reflection as
an algebra automorphism, with the spherical line $R_0$ fixed.  Therefore
the products
\[
 u_+^{\,k-j}u_-^{\,j}\in R_k,\qquad 0\leq j\leq k,
\]
are highest vectors of weight $k\varpi_2$.  They are nonzero because
$R$ is a domain, and linearly independent because their circle weights
$2k-4j$ are distinct.  The dimension calculation above shows that they
form a basis of the highest-vector space.  Finally,
\[
 r\bigl(u_+^{\,k-j}u_-^{\,j}\bigr)=u_+^{\,j}u_-^{\,k-j}.
\]
The monomial correspondence $z_+^{k-j}z_-^j\mapsto u_+^{k-j}u_-^j$
therefore gives an $O(2)$-isomorphism
\[
 \Sym^k Z\simeq\Hom_{K_F}(T_k,R_k).
\]
For $k=0$ this is the fixed spherical line; for $k=2a>0$ the
zero-weight line is spanned by $(u_+u_-)^a$ and is fixed by $r$.
Together with freeness, this proves the multiplication isomorphism
in part~(1).

For completeness, the calculation for $Q_F$ is similar but simpler.  The
$\mathrm{SL}(W)$-type $W$ in \cref{eq:Cauchy-rectangle} occurs when
$s=r-1$, $1\leq r\leq n$.  The second assertion of
\cref{lem:paired-branching} shows that each corresponding factor
contains $V_{\mathrm{Sp}(3)}(\varpi_1)$ once on restriction.
Thus $\dim\Hom_{K_F}(Q_F,R_n)=n$ for $n\geq1$, and
\[
  H_{M_{Q_F}}(t)=
 \frac{t}
 {(1-t)^2}.
\]
 Finally,
\cref{eq:70} says that its degree-one generator is the sign, which supplies
the factor $\det$ in part (2).
\end{proof}

\pagebreak[3]
\begin{theorem}[Compact-type multiplicity spaces]\label{thm:multiplicity}
The following hold.
\begin{enumerate}
\item For $T_k=V_{\mathrm{Sp}(3)}(k\varpi_2)\boxt\mathbf1$,
\begin{equation}\label{eq:Mk-induced}
 \mathcal M(T_k)\simeq
 \Ind_{O(2)}^{(\s,O(2))}\bigl(\Sym^kX\bigr).
\end{equation}
\item For $Q_F=V(\varpi_1)\boxt W$,
\[
\mathcal M(Q_F)\simeq\Ind_{O(2)}^{(\s,O(2))}\det.
\]
\end{enumerate}
\end{theorem}

\begin{proof}
Choose the lowest-degree $O(2)$-stable generator space supplied by
\cref{lem:covariant-branching}, and lift it to the same filtration degree
in $\mathcal M(\kappa)$.  The enveloping-algebra action gives a filtered map
\[
 U(\s)\otimes_{U(\mathfrak{so}_2)}E\longrightarrow\mathcal M(\kappa).
\]
By PBW its associated graded is the multiplication isomorphism
$\Sym(X)\otimes E\simeq M_\kappa$ of that lemma, since $X\simeq Z$.
The filtrations are bounded below and exhaustive, so the filtered map is
an isomorphism.  Take $E=\Sym^kX$ or $E=\det$, respectively.
\end{proof}

Note that the $O(2)$ generator in part (1) is
\begin{equation}\label{eq:Ek}
E_k:=\Sym^kX,\qquad E_k|_{SO(2)}=\bigoplus_{j=0}^k\C_{2k-4j}.
\end{equation}
Pairing opposite nonzero weights and retaining the reflection on the zero
line gives
\begin{equation}\label{eq:Ek-O2}
E_{2a}=\mathbf1\oplus\tau_4\oplus\tau_8\oplus\cdots\oplus\tau_{4a},
\quad
E_{2a+1}=\tau_2\oplus\tau_6\oplus\cdots\oplus\tau_{4a+2}.
\end{equation}
In particular $E_2$ has weights $-4,0,4$, and its zero-weight line is the
\emph{trivial}, not the sign, representation of $O(2)$.
More generally, the zero-weight line in every $E_{2a}$ is trivial,
by the product calculation in the proof of \cref{lem:covariant-branching}.

\section{See--saw reciprocity and ping--pong}
\label{sec:pingpong}

We combine the compact-type lifts computed in the preceding two sections
using see--saw reciprocity.  The resulting ping--pong argument proves
uniqueness of the irreducible theta quotient, with the spherical case
completed by the exceptional-fibre analysis in
\cref{sec:exceptional-fibre}.  We also establish the Capelli identity
relating the two Casimir operators, which is needed for that analysis.

\subsection{Big theta lifts and adjunction}

If $G_1$ is one member of a dual pair and $Y$ is an admissible $G_1$-module, set
\[
 \ThetaB(Y)=(\PiMin\otimes Y^\vee)_{G_1}.
\]
The maximal semisimple quotient is denoted $\thetaS(Y)$.  For a compact type
we use the same notation for its isotypic multiplicity space.  Algebraic
Frobenius reciprocity, first for $H\times F$ and then for the commuting pair
$O(2)\times F$, gives
\begin{align}
 \Hom_{K_F}\bigl(\ThetaB_F(\sigma),\kappa\bigr)
  &\simeq
 \Hom_{(\s,O(2))}\bigl(\mathcal M(\kappa),\sigma\bigr),
                                                        \label{eq:ss1}\\
 \Hom_{O(2)}\bigl(\ThetaB_H(\pi),\tau\bigr)
  &\simeq
 \Hom_{(\f,K_F)}\bigl(\mathcal A_\tau,\pi\bigr).
                                                        \label{eq:ss2}
\end{align}
The second identity uses \cref{eq:O2-F-module}.
By \cref{thm:multiplicity} and Frobenius reciprocity once more,
\begin{equation}\label{eq:type-identity}
 \Hom_{K_F}\bigl(\ThetaB_F(\sigma),T_k\bigr)
       \simeq\Hom_{O(2)}(E_k,\sigma).
\end{equation}
The analogous identities for $\mathbf1$ and $Q_F$ use respectively
$\mathbf1$ and $\det$.

\subsection{The ping--pong lemma}

We isolate the elementary mechanism.  It is useful well beyond the present
dual pair.

\begin{lemma}[Ping--pong]\label[lemma]{lem:pingpong}
Let $\sigma$ be irreducible, let $\tau$ be an $O(2)$-type occurring
in $\sigma$, and let $\mathcal F_\tau$ generate $\mathcal A_\tau$.
Suppose that $\ThetaB_F(\sigma)$ has finite length, that a $K_F$-type
$\kappa$ occurs in it with multiplicity one, and that no other constituent
of $\mathcal F_\tau$ occurs in it.  Then its irreducible quotient is unique
and contains $\kappa$ with multiplicity one.
\end{lemma}

\begin{proof}
The multiplicity assumption makes the big lift nonzero, so finite length
gives an irreducible quotient $L$.  Reciprocity applied to this quotient
gives a nonzero map $\ThetaB_H(L)\to\sigma$, which is surjective.
Since $\tau$ occurs in $\sigma$, exactness of compact isotypic projection
and \cref{eq:ss2} give a nonzero map $\mathcal A_\tau\to L$.
It is surjective because $L$ is irreducible.  Consequently $L$ meets
$\mathcal F_\tau$, and the support assumption forces it to contain
$\kappa$.  Every simple summand of the cosocle has this property.
The multiplicity-one hypothesis therefore makes that cosocle simple.
\end{proof}

\begin{proof}[Proof of \cref{thm:main}]
If $\sigma$ is $O(2)$-spherical, \cref{thm:unique-exception} shows that every
irreducible quotient of its big lift is spherical unless
$\sigma\simeq S^4(\C^2)$.  In the nonexceptional case the trivial $K_F$-type
has multiplicity one by \cref{thm:multiplicity}(1), so the quotient is
unique.  At the exceptional parameter, \cref{thm:exceptional-fibre} gives
exactly the two quotients $\Pi_6^+$ and $\Pi_6^-$.  This proves (1).
Part (2) uses $Q_F$ and the sign line in
\cref{thm:multiplicity}(2).  It too gives uniqueness among all irreducible
quotients.

Suppose next that the minimal $O(2)$-type is $\tau_{2k}$, with $k\geq1$.
Take $\tau=\tau_{2k}$ and
$\kappa=T_k=V(k\varpi_2)\boxt\mathbf1$.  Formula \cref{eq:Ek-O2} shows
that $\tau_{2k}$ first occurs in $E_k$, with multiplicity one.  On the
forward side,
\[
                 \mathcal F_k=T_0\oplus T_1\oplus\cdots\oplus T_k.
\]
For $j<k$, every $O(2)$-type in $E_j$ has first occurrence strictly below
$k$ and is therefore absent from $\sigma$.  Thus all the hypotheses of
\cref{lem:pingpong} hold, proving (3).

Finally assume that $\sigma$ is obtained by adjoining the $O(2)$ reflection
to the lowest- and highest-weight modules with lowest weight $2k$, where
$k\geq1$. Theorem \ref{thm:lowest-weight-lift} identifies its unique quotient with
$A_{\mathfrak q}((k-6)e_1)$; equivalently, one may use the opposite
parabolic for the highest-weight orientation.  This proves (4).
\end{proof}

\subsection{Capelli identity}
\label{sec:capelli}

Let $\Omega_H$ be the Casimir of $\s$, normalized to act by $x^2-1$ at
Harish--Chandra parameter $x$, and let $\Omega_F$ be the Casimir of $\f$,
normalized by
\[
       \Omega_F(\lambda)=\langle\lambda,\lambda\rangle
                         -\langle\rho_{F_4},\rho_{F_4}\rangle.
\]
The infinitesimal-character transfer for this dual pair is
\begin{equation}\label{eq:ic-transfer}
             x\longmapsto\tfrac12(x,5,3,1);
\end{equation}
see \cite{LiIC,SavinPingPong}.  It gives
\[
 \Omega_F\left(\tfrac12(x,5,3,1)\right)=\frac{x^2-121}{4}.
\]
In fact this gives the following identity of operators on $\PiMin$:
\begin{equation}\label{eq:casimir-relation}
                    \Omega_H=4\Omega_F+120.
\end{equation}

To justify the operator identity, let $J_{\mathrm{Jos}}$ be the Joseph ideal in
$U(\mathfrak e_{7,\C})$.  The action on any minimal representation
factors through $U(\mathfrak e_{7,\C})/J_{\mathrm{Jos}}$; since $J_{\mathrm{Jos}}$ is maximal,
this quotient acts faithfully on any such nonzero representation.
We may therefore check the identity on a minimal representation for
the real form in which the $F_4$ member is compact.  There the known
decomposition under the dual pair and its infinitesimal characters
\cite{LiIC} show that $\Omega_H-4\Omega_F-120$ acts by zero on every
summand.  Thus this element belongs to $J_{\mathrm{Jos}}$.
The dual pair $\mathfrak{sl}_2\oplus\mathfrak f_4$ is unique up to
conjugacy over $\C$.  This conjugacy identifies the normalized Casimirs
and preserves the Joseph ideal, so the same identity holds on the
split minimal representation $\PiMin$.

\section{The exceptional spherical fibre and its uniqueness}
\label{sec:exceptional-fibre}

We show that the spherical five-dimensional representation $S^4(\C^2)$
has two irreducible theta quotients, obtained by restricting the minimal
representation of split $E_6$ to $F_{4(4)}$.  This gives an archimedean
counterpart of the exceptional two-valued fibre in
\cite{KarasiewiczSavin}.  We then use the Capelli identity and a
calculation in the Pl\"ucker ring to prove that every other spherical
parameter has a unique irreducible theta quotient.

\subsection{The fibre at \texorpdfstring{$s=5/2$}{s=5/2}}
\label{sec:exceptional-value}

Let $G_6^{\mathrm{sp}}=E_{6(6)}$ be the split group of type $E_6$.  This is
the derived group of the relevant Levi subgroup of split $E_7$; it should
not be confused with the quasi-split group $E_{6(2)}$ occurring in the
see--saw.  Write $\Pi_6$ for the spherical minimal representation of
$G_6^{\mathrm{sp}}$.  The nontrivial diagram involution $\iota$ of
$G_6^{\mathrm{sp}}$ has fixed group $F_{4(4)}$.  After extending $\Pi_6$ to
$G_6^{\mathrm{sp}}\rtimes\langle\iota\rangle$ so that $\iota$ fixes the
spherical vector, write
\begin{equation}\label{eq:E6-eigenspaces}
 \Pi_6\big|_{\langle\iota\rangle\times F_{4(4)}}
   =\bigl(\mathbf1\boxt\Pi_6^+\bigr)
       \oplus\bigl(\operatorname{sgn}\boxt\Pi_6^-\bigr).
\end{equation}
The two $F_{4(4)}$-modules in \cref{eq:E6-eigenspaces} are irreducible and
inequivalent by \cref{thm:E6-F4-restriction}.  Its proof in
Section~\ref{sec:E6-F4-restriction} uses the Lagrangian Grassmannian
linear section and unitarity, independently of the present argument.
For the analogous split $p$-adic restriction, see
\cite{KarasiewiczSavin}.

Their minimal $K_F=\mathrm{Sp}(3)\times\mathrm{SU}(2)$-types can be read
off in the first two $K_{G_6^{\mathrm{sp}}}=\mathrm{Sp}(4)$-types of
$\Pi_6$.  The degree-zero type is the trivial representation, while
\begin{equation}\label{eq:Sp4-first-branching}
 \left.V_{\mathrm{Sp}(4)}(\varpi_4)\right|_{\mathrm{Sp}(3)\times
 \mathrm{SU}(2)}
 =\bigl(V_{\mathrm{Sp}(3)}(\varpi_3)\boxt\C^2\bigr)
       \oplus
   \bigl(V_{\mathrm{Sp}(3)}(\varpi_2)\boxt\mathbf1\bigr).
\end{equation}
The first summand in \cref{eq:Sp4-first-branching} is the tangent
representation $\mathfrak p_F$ of $F_{4(4)}/K_F$ and is fixed by $\iota$;
the second is the normal representation and is negated by $\iota$.  Thus
\begin{equation}\label{eq:Pi6-minimal-types}
 \operatorname{minK}(\Pi_6^+)=T_0=\mathbf1,
 \qquad
 \operatorname{minK}(\Pi_6^-)=T_1
      =V_{\mathrm{Sp}(3)}(\varpi_2)\boxt\mathbf1.
\end{equation}

We will also need that $Q_F=V_{\mathrm{Sp}(3)}(\varpi_1)\boxt\C^2$
does not occur anywhere in $\Pi_6$.  The compact branching formula
\eqref{eq:E6-compact-branching} describes all its $K_F$-types.
Every $\mathrm{Sp}(3)$ weight in that formula has zero $\varpi_1$-coefficient,
so $Q_F$ is absent from both $\Pi_6^+$ and $\Pi_6^-$.

Let $B_H$ be a Borel subgroup of $H=\mathrm{PGL}_2(\R)$, and normalize the
spherical principal series by
\[
                 I(s)=i_{B_H}^{H}(|\cdot|^s).
\]
Our convention is the one for which the positive reducibility point
$s=5/2$ has finite-dimensional quotient $S^4(\C^2)$.  If $D_6$ denotes the
irreducible $(\s,O(2))$-module obtained by adjoining the reflection to the
holomorphic and antiholomorphic discrete series of lowest $K$-type $6$,
then
\begin{equation}\label{eq:PGL2-reducibility}
 0\longrightarrow D_6\longrightarrow I(5/2)
   \longrightarrow S^4(\C^2)\longrightarrow0.
\end{equation}
In particular,
\begin{equation}\label{eq:PGL2-types-exceptional}
 \left.S^4(\C^2)\right|_{O(2)}
       =\tau_0^+\oplus\tau_2\oplus\tau_4,
 \qquad
 \left.D_6\right|_{O(2)}=\tau_6\oplus\tau_8\oplus\tau_{10}\oplus\cdots.
\end{equation}

\begin{theorem}[The exceptional fibre]\label{thm:exceptional-fibre}
For the spherical finite-dimensional representation $S^4(\C^2)$ of
$\mathrm{PGL}_2(\R)$, the small theta lift is
\begin{equation}\label{eq:exceptional-small-theta}
       \thetaS_F\bigl(S^4(\C^2)\bigr)
             \simeq \Pi_6^+\oplus\Pi_6^-.
\end{equation}
Consequently its big theta lift has exactly two irreducible quotients.  The
first is spherical, and the second has minimal $K_F$-type
$V_{\mathrm{Sp}(3)}(\varpi_2)\boxt\mathbf1$.  In particular the theta
correspondence is not one-to-one at this parameter.
\end{theorem}

\begin{proof}
Let $P_7=M_7N_7$ be the maximal parabolic subgroup of split $E_7$ whose
derived Levi subgroup is $G_6^{\mathrm{sp}}$.  The global minimal
representation of $E_7$ is the automorphic theta representation constructed
by Ginzburg--Rallis--Soudry.  Their constant-term calculation shows that its
constant term along $P_7$ supplies the global minimal representation
of the $E_6$ factor and a character of the Levi
\cite{GinzburgRallisSoudry}.  Passing to the real component gives,
with unnormalized coinvariants, a surjection
\begin{equation}\label{eq:E7-E6-Jacquet}
       (\Pi_7)_{N_7}
          \twoheadrightarrow \Pi_6\otimes|\cdot|^3\ \oplus\ |\cdot|^6.
\end{equation}
Here $\Pi_7=\PiMin$; only projection to the first summand is used.
Surjectivity is sufficient for the following Frobenius-reciprocity
argument; the full Jacquet module need not be determined.

The split torus in the centre of $M_7$ is the split Cartan of the rank-one
member $H$.  Projection to either $F_{4(4)}$-summand and Frobenius
reciprocity give nonzero induced-representation maps.  Initially allow
the two extensions of the connected principal series: write
$I_\varepsilon(5/2)=I(5/2)\otimes\eta_\varepsilon$, where
$\eta_\varepsilon$ is either $\mathbf1$ or the character $\det$ of
$H/H^\circ$.  Thus we have nonzero $H\times F_{4(4)}$-maps
\begin{equation}\label{eq:Frobenius-exceptional}
              \Pi_7\longrightarrow I_\varepsilon(5/2)\boxt\Pi_6^\varepsilon,
              \qquad \varepsilon\in\{+,-\}.
\end{equation}
We will show that both $\eta_\varepsilon$ are trivial.

Neither summand can be attached to $D_6$.  Indeed,
\cref{eq:type-identity,eq:Ek-O2,eq:PGL2-types-exceptional} give
\begin{align*}
 \Hom_{K_F}\bigl(\ThetaB_F(D_6),T_0\bigr)
    &\simeq\Hom_{O(2)}(\tau_0^+,D_6)=0,\\
 \Hom_{K_F}\bigl(\ThetaB_F(D_6),T_1\bigr)
    &\simeq\Hom_{O(2)}(\tau_2,D_6)=0.
\end{align*}
By \cref{eq:Pi6-minimal-types}, neither $\Pi_6^+$ nor $\Pi_6^-$ can
therefore be a quotient of $\ThetaB_F(D_6)$.  Twisting $D_6$ by $\det$
does not change its isomorphism class (or its displayed compact types).
Consequently composition in \cref{eq:Frobenius-exceptional} with
\[
 I_\varepsilon(5/2)\twoheadrightarrow
           S^4(\C^2)\otimes\eta_\varepsilon
\]
is nonzero: otherwise the image would lie in the discrete-series
submodule, contrary to the preceding type calculation.

It remains to exclude the sign-spherical five-dimensional representation
$S^4(\C^2)\otimes\det$.  If either $\Pi_6^\varepsilon$ were its theta
quotient, reverse reciprocity \cref{eq:ss2} would give a nonzero, hence
surjective, map
\[
             \mathcal A_{\tau_0^-}\longrightarrow\Pi_6^\varepsilon.
\]
By \cref{thm:finite-generation}, the source is generated by $Q_F$.
Its image must therefore contain $Q_F$, contradicting the absence of
that type from all of $\Pi_6$.  Hence $\eta_+=\eta_-=\mathbf1$,
and both summands are theta quotients of the spherical $S^4(\C^2)$.
For $\Pi_6^+$ this also follows from its $K_F$-fixed vector, namely the
spherical vector of $\Pi_6$: a sign-spherical parameter has no trivial
$K_F$-type in its lift by \cref{eq:type-identity} with $k=0$.
This produces the two summands on the right side of
\cref{eq:exceptional-small-theta}.

It remains only to exclude a third quotient.  By
\cref{thm:finite-generation}, every irreducible quotient over an
$O(2)$-spherical parameter meets the finite generator $T_0\oplus T_1$.
For the parameter at hand, \cref{eq:type-identity} gives
\begin{align*}
 \dim\Hom_{K_F}\bigl(\ThetaB_F(S^4(\C^2)),T_0\bigr)
   &=\dim\Hom_{O(2)}\bigl(\tau_0^+,S^4(\C^2)\bigr)=1,\\
 \dim\Hom_{K_F}\bigl(\ThetaB_F(S^4(\C^2)),T_1\bigr)
   &=\dim\Hom_{O(2)}\bigl(\tau_2,S^4(\C^2)\bigr)=1.
\end{align*}
Hence at most one quotient can contain either generator.  The two quotients
already constructed contain $T_0$ and $T_1$, respectively, so they exhaust
the cosocle.
\end{proof}

\begin{remark}
The same exceptional value $s=5/2$ appears in the split $p$-adic
correspondence.  There the small theta lift is the two-summand cosocle of an
$F_4$ degenerate principal series; see \cite{KarasiewiczSavin}.  The argument
above identifies the real analogues of those two summands by their minimal
$K_F$-types.
\end{remark}

\subsection{Uniqueness of the exceptional infinitesimal character}
\label{sec:unique-exception}

We now prove that the fibre in \cref{thm:exceptional-fibre} is the only
failure of one-to-one correspondence over an $O(2)$-spherical parameter.
Put
\begin{equation}\label{eq:A-and-B}
 \mathcal A=\mathcal A_{\tau_0^+}=\Hom_{O(2)}(\mathbf1,\PiMin),
 \qquad
 \mathcal B=U(\f)v_0\subset \mathcal A,
\end{equation}
where $v_0$ is the spherical vector.  Thus $\mathcal A$ is generated by
\[
 T_0=\mathbf1,
 \qquad
 T_1=V_{\mathrm{Sp}(3)}(\varpi_2)\boxt\mathbf1,
\]
whereas $\mathcal B$ is the submodule generated by $T_0$.  For a $K_F$-type $T$ we
write
\[
             \mathcal A[T]=\Hom_{K_F}(T,\mathcal A),
             \qquad \mathcal B[T]=\Hom_{K_F}(T,\mathcal B)
\]
for the corresponding multiplicity spaces.  Equip $\mathcal A$ with the
filtration inherited from $\PiMin$, and $\mathcal B$ with the induced filtration;
the multiplicity spaces carry the corresponding filtrations.

The failure of one-to-one correspondence in
\cref{thm:exceptional-fibre} already shows that
$\mathcal B[T_1]\subsetneq \mathcal A[T_1]$.  Indeed, the nonspherical quotient
$\Pi_6^-$ gives a nonzero map $\mathcal A\to\Pi_6^-$ by see--saw reciprocity.
It kills $v_0$ and hence $\mathcal B$, since $\Pi_6^-$ has no trivial $K_F$-type.
It is nonzero on the generating type $T_1$.  Our goal is to prove
\[
                   \dim\bigl(\mathcal A[T_1]/\mathcal B[T_1]\bigr)=1.
\]
It suffices to prove this for the associated graded spaces: the induced
filtration on $\mathcal B[T_1]$ gives
\[
 \gr\bigl(\mathcal A[T_1]/\mathcal B[T_1]\bigr)
       =\gr \mathcal A[T_1]/\gr \mathcal B[T_1].
\]

Put $z=z_+z_-$.  Taking $O(2)$-invariants in
Lemma~\ref{lem:covariant-branching}(1), with $k=1$, shows that
\begin{equation}\label{eq:grAT1-free}
                   \gr \mathcal A[T_1]=\C[z]\,\bar t_1,
                   \qquad \deg\bar t_1=2,\quad\deg z=2.
\end{equation}
In particular, this is a free rank-one $\C[z]$-module, with one line
in each even degree beginning in degree two.  The principal symbol of
the $\mathfrak{sl}_2$ Casimir $\Omega_H$ is a nonzero multiple of $z$.
The Capelli identity \cref{eq:casimir-relation} implies that $\mathcal B$ is
stable under $\Omega_H$, because it is a $U(\f)$-submodule.
Thus $\gr \mathcal B[T_1]$ is a graded $\C[z]$-submodule of
\cref{eq:grAT1-free}.  Once it is nonzero, it has the form
\begin{equation}\label{eq:grBT1-ideal}
                   \gr \mathcal B[T_1]=z^d\C[z]\,\bar t_1
\end{equation}
for an integer $d\geq1$; properness excludes $d=0$.
The desired dimension is precisely $d$.  The following lemma proves
nonvanishing and $d=1$.

\begin{lemma}[The degree-four bridge]\label[lemma]{lem:degree-four-bridge}
Let $P=V(\varpi_3)\boxt W=(\mathfrak p_F)_\C\subset R_1$, and let
\[
                     m_d:\Sym^dP\longrightarrow R_d
\]
be Cartan multiplication in the Pl\"ucker ring.  Then
\[
 \Hom_{K_F}(T_1,\operatorname{im}m_2)=0,
 \qquad
 \dim\Hom_{K_F}(T_1,\operatorname{im}m_4)=1.
\]
Consequently
\[
             \gr \mathcal B[T_1]=z\C[z]\,\bar t_1,
             \qquad
             \dim\bigl(\gr \mathcal A[T_1]/\gr \mathcal B[T_1]\bigr)=1.
\]
\end{lemma}

\begin{proof}
The reflection fixes $P$ pointwise by \cref{eq:70}, so every
$\operatorname{im}m_d$ is $O(2)$-fixed.  Put $E=\Lambda^3_0V$, so
$P=E\boxt W$.  The branching identity
\[
       (\Sym^2P)^{\mathrm{SU}(2)}=\Lambda^2E
          =\mathbf1\oplus V(2\varpi_2)
\]
shows directly that $T_1$ does not occur in $\operatorname{im}m_2$.
We construct the degree-four copy by squaring the $V(2\varpi_2)$
summand, using the $AC$-coordinates.

\emph{The quadratic copy of $V(2\varpi_2)$.}
In the following calculation, $A,B,C$ denote the degree-one summands
of \cref{eq:ABC-decomposition}, and products are taken in $R$.
The Littlewood--Richardson rule gives
\[
 (R_2)^{\mathrm{SU}(2)}_{\text{$SO(2)$-weight }0}
       \simeq V_{\mathrm{SU}(6)}(2,2,1,1,0,0),
\]
which is isomorphic to the $189$-dimensional Cartan summand of
$A\otimes C$.  The product of highest weight vectors of the
irreducible $\mathrm{SU}(V)$-modules $A$ and $C$ is nonzero, since
$R$ is a domain, and generates this summand.  Together with the
inclusion $AC\subset(B^2)^{\mathrm{SU}(2)}$ from
\cref{prop:B-linear-section}, this proves
\[
       AC=(B^2)^{\mathrm{SU}(2)}
          \simeq V_{\mathrm{SU}(6)}(2,2,1,1,0,0).
\]
Its restriction to $\mathrm{Sp}(3)$ is
\[
 \mathbf1\oplus2V(\varpi_2)\oplus V(2\varpi_2)
                   \oplus V(\varpi_1+\varpi_3).
\]
In particular, $AC$ contains a unique copy of $V(2\varpi_2)$.
A highest weight vector of this type has the two expressions
\begin{align*}
 p_{e_1e_2e_3f_3}\,p_{e_1e_2xy}
   ={}&p_{e_1e_2e_3x}\,p_{e_1e_2f_3y}
       -p_{e_1e_2e_3y}\,p_{e_1e_2f_3x}.
\end{align*}
Use the Pl\"ucker relation with triple $(e_1,e_2,x)$ and
quintuple $(e_1,e_2,e_3,f_3,y)$.  The left side belongs to $AC$;
the right side belongs to $(\Sym^2P)^{\mathrm{SU}(2)}$, since
$e_1\wedge e_2\wedge e_3$ and $e_1\wedge e_2\wedge f_3$ are
primitive three-forms.  The left side is nonzero because $R$ is a
domain.  Thus the copy of $V(2\varpi_2)$ in $AC$ lies in
$\operatorname{im}m_2$.

\emph{The Jordan algebra and the projection onto $V(2\varpi_2)$.}
Let $\mathcal J$ be the Jordan algebra of symplectically self-adjoint
endomorphisms of $V$, with
\[
 X\circ Y=\tfrac12(XY+YX),\qquad
 \tau(X)=\operatorname{tr}_V(X),\qquad
 (X,Y)=\tau(X\circ Y).
\]
There is an equivariant identification $\Lambda^2V\simeq\mathcal J$
given by
\[
 u\wedge v\longmapsto
       \bigl(z\longmapsto\langle z,v\rangle u-\langle z,u\rangle v\bigr).
\]
It sends $\Omega_V=\sum_i e_i\wedge f_i$ to $I_V$ and has trace
$2\langle u,v\rangle$ on $u\wedge v$.  Hence the primitive
submodule $\Lambda^2_0V$ identifies with
\[
 T=V(\varpi_2)=\{X\in\mathcal J:\tau(X)=0\},\qquad
                   \mathcal J=\C I_V\oplus T.
\]
Write $X_0=X-\tfrac16\tau(X)I_V$ for the traceless part.
Jordan multiplication induces
\[
 m:\Sym^2T\longrightarrow\mathcal J,\qquad m(X\odot Y)=X\circ Y.
\]
Since $\Sym^2T=\mathbf1\oplus T\oplus V(2\varpi_2)$ and $m$ is
nonzero on the first two summands, we have
\[
                         \ker m=V(2\varpi_2).
\]
Thus the projection $\pi_{90}$ simply removes the scalar and
$T$ components.  More explicitly, equip $\Sym^2T$ with the pairing
\[
 (X\odot Y,U\odot Z)_2
       =\tfrac12\bigl((X,U)(Y,Z)+(X,Z)(Y,U)\bigr),
\]
and let $m^*$ be the adjoint of $m$ for these trace pairings.
For trace-dual bases $t_i,t^i$ of $T$, it is given by
\[
                    m^*(Z)=\sum_i t_i\odot(Z\circ t^i)_0.
\]
A calculation using Jordan multiplication gives
\[
 mm^*|_{\C I_V}=\tfrac73\operatorname{id},\qquad
 mm^*|_T=\tfrac23\operatorname{id},\qquad
 \pi_{90}=\operatorname{id}-m^*(mm^*)^{-1}m.
\]

\emph{The quartic by trace contraction.}
For $u,v,s,t\in T$, define a symmetric bilinear operation on
$\Sym^2T$ by
\begin{align*}
 (u\odot v)\star(s\odot t)=\tfrac14\bigl(&
 (v,s)u\odot t+(v,t)u\odot s\\
 &+(u,s)v\odot t+(u,t)v\odot s\bigr).
\end{align*}
This contracts one factor from each tensor using the trace pairing.
It gives a basis-independent meaning to squaring a quadratic tensor,
without regarding it as an endomorphism of $T$.

For a four-vector, take its $A$-component in $\Lambda^4V$ and its
$C$-component in $\Lambda^2V\otimes\Lambda^2W$.  Symplectically
contract the former and identify the latter with $\Lambda^2V$ using
$x\wedge y$.  Their images in $\mathcal J$ need not be traceless;
let $a,c\in T$ be their traceless parts.  Set
\begin{equation}\label{eq:bridge-quartic}
            D=\pi_{90}(a\odot c),\qquad
            \Psi=\bigl(m(D\star D)\bigr)_0.
\end{equation}
The coordinate functions of $D$ belong to the unique
$V(2\varpi_2)$ in $AC$, which the Pl\"ucker identity above identifies
with the image of $V(2\varpi_2)\subset\Sym^2P$.
Thus the coordinates of $\Psi$ belong to $\operatorname{im}m_4$.
The construction is $\mathrm{Sp}(3)$-equivariant and
$\mathrm{SU}(2)$-invariant.

\emph{A nonzero value on the Pl\"ucker cone.}
Let $E_i\in\mathcal J$ be the image of $e_i\wedge f_i$.  These are
orthogonal idempotents: $E_i\circ E_j=\delta_{ij}E_i$,
$\sum_iE_i=I_V$, and $\tau(E_i)=2$.
At the decomposable vector
\[
           \xi=(e_1+x)\wedge(f_1+y)\wedge e_2\wedge f_2,
\]
the two components in $\mathcal J$ are $E_1+E_2$ and $E_2$.
Consequently their traceless parts are
\[
       a=\tfrac13(E_1+E_2-2E_3),\qquad
       c=\tfrac13(-E_1+2E_2-E_3).
\]
Using the formulas for $\pi_{90}$ and $\star$ above gives
\begin{equation}\label{eq:bridge-evaluation}
                  \Psi(\xi)=\tfrac1{16}(-2E_1+E_2+E_3)\ne0.
\end{equation}
The wedge pairing identifies $\Lambda^4U$ with its dual and is
nondegenerate on $P$.  Thus this evaluation proves that the
coordinate functions of $\Psi$ span a nonzero copy of $T_1$ in
$\operatorname{im}m_4$.  By \cref{lem:covariant-branching}(1),
$T_1$ has multiplicity one in $R_4^{O(2)}$.  This proves the
second assertion.

Return now to the filtered module $\mathcal B=U(\f)v_0$.
Lift this nonzero copy of $T_1$ in $\operatorname{im}m_4$ by PBW
symmetrization and apply it to $v_0$.  Its leading symbol is the chosen
copy in $R_4$, so $\gr \mathcal B[T_1]$ contains the degree-four line
$\C z\bar t_1$.  Therefore the exponent in \cref{eq:grBT1-ideal}
satisfies $d\leq1$.  Since $\mathcal B[T_1]$ is proper in $\mathcal A[T_1]$, we have
$d\geq1$, and hence $d=1$.  This argument uses only the leading symbol
of the chosen operator; it does not identify $\gr \mathcal B$, for its induced
filtration, with the image of $\Sym P\to R$.
\end{proof}

Returning to the filtered multiplicity spaces, choose a lift $t_1$ of
$\bar t_1$.  The symbol calculation in \cref{eq:grAT1-free} gives
\begin{equation}\label{eq:AT1-free}
                       \mathcal A[T_1]=\C[\Omega_H]t_1.
\end{equation}
Indeed, the vectors $\Omega_H^j t_1$ have nonzero leading symbols in
the successive one-dimensional graded pieces.  They are therefore
linearly independent and span by induction on filtration degree.
By \cref{lem:degree-four-bridge}, $\mathcal B[T_1]$ has codimension one in
this free module.  Since it is stable under $\Omega_H$, there is a
scalar $c$ such that
\begin{equation}\label{eq:BT1-polynomial}
                  \mathcal B[T_1]=(\Omega_H-c)\C[\Omega_H]t_1.
\end{equation}

The map $\mathcal A\to\Pi_6^-$ used above is nonzero on $\mathcal A[T_1]/\mathcal B[T_1]$.
It comes from the theta quotient of $S^4(\C^2)$, whose
Harish--Chandra parameter is $x=5$.  Thus $\Omega_H$ acts on this
one-dimensional quotient by $5^2-1=24$, and $c=24$.  The Capelli
identity now gives
\begin{equation}\label{eq:p-explicit}
 \mathcal B[T_1]=(\Omega_H-24)\mathcal A[T_1]=(\Omega_F+24)\mathcal A[T_1],
 \qquad \Omega_H-24=4(\Omega_F+24).
\end{equation}

\begin{proposition}[Casimir localization]\label[proposition]{prop:casimir-localization}
The quotient $\mathcal A/\mathcal B$ is annihilated by $\Omega_F+24$.  In particular,
\begin{equation}\label{eq:localized-generation}
        \mathcal A[(\Omega_F+24)^{-1}]=\mathcal B[(\Omega_F+24)^{-1}].
\end{equation}
Thus, away from the infinitesimal character $\Omega_F=-24$, the trivial
$O(2)$-multiplicity module $\mathcal A$ is generated by its spherical vector alone.
\end{proposition}

\begin{proof}
Theorem \ref{thm:finite-generation} says that $\mathcal A/\mathcal B$ is generated by the
image of $T_1$.  Equation \eqref{eq:p-explicit} shows that $\Omega_F+24$ annihilates that image, and
$\Omega_F$ is central in $U(\f)$.  It therefore annihilates all of $\mathcal A/\mathcal B$.
\end{proof}

\begin{theorem}[Uniqueness of the exceptional fibre]
\label{thm:unique-exception}
Let $\sigma$ be an irreducible $(\s,O(2))$-module containing
$\tau_0^+$.  If $\sigma\not\simeq S^4(\C^2)$, then
$\ThetaB_F(\sigma)$ has exactly one irreducible quotient and it is
spherical.  For $\sigma\simeq S^4(\C^2)$ it has exactly the two quotients
$\Pi_6^+$ and $\Pi_6^-$ of \cref{thm:exceptional-fibre}.
\end{theorem}

\begin{proof}
Let $x$ be the Harish--Chandra parameter of $\sigma$, and let $L$ be an
irreducible quotient of $\ThetaB_F(\sigma)$.  Theta adjunction gives a
nonzero map $\ThetaB_H(L)\to\sigma$.  It is surjective, and exactness of
taking $O(2)$-types shows that the trivial type in $\sigma$ is already
present in $\ThetaB_H(L)$.  Now \cref{eq:ss2}, with $\tau=\tau_0^+$, gives
a nonzero $F_{4(4)}$-map $\mathcal A\to L$.  On $L$ the central element in
\cref{eq:localized-generation} acts by
\[
                       \Omega_F+24=\frac{x^2-25}{4}.
\]
If $x^2\ne25$, \cref{prop:casimir-localization} implies that the map
$\mathcal A\to L$ is already generated by the image of $v_0$.  Hence $L$ contains
$T_0$ and is spherical.  The $T_0$-multiplicity in the big lift is one by
\cref{thm:multiplicity}(1).  Since every irreducible summand of its cosocle
contains $T_0$, that cosocle has length one.  Thus the irreducible quotient
is unique.

It remains to consider $x^2=25$.  At this infinitesimal character the two
irreducible subquotients of the spherical principal series are
$S^4(\C^2)$ and $D_6$ as in \cref{eq:PGL2-reducibility}; the latter contains
no $\tau_0^+$.  Thus an irreducible $O(2)$-spherical $\sigma$ must be
$S^4(\C^2)$, and \cref{thm:exceptional-fibre} gives precisely the two stated
quotients.
\end{proof}

\section{Discrete-series lifts via the unitary correspondence}
\label{sec:lowest-weight}

We deduce part~(4) of \cref{thm:main} from the discrete summands in
\cite[Theorems~A and~B(2)]{FrahmLabriet}.  Write $\sigma_{2k}$ for the
irreducible $(\s,O(2))$-module whose restriction to $(\s,SO(2))$ is the sum
of the lowest- and highest-weight modules with extremal weights
$2k$ and $-2k$, where $k\geq1$.

We first specify the parabolic and compare normalizations.  Realize
the $F_4$ roots in an orthonormal basis $e_1,e_2,e_3,e_4$, and put
\[
 \epsilon_1=\tfrac12(e_1+e_2),\qquad
 \epsilon_2=\tfrac12(e_1-e_2),\qquad
 \epsilon_3=\tfrac12(e_3+e_4),\qquad
 \delta=\tfrac12(e_3-e_4).
\]
The compact roots are $\pm2\epsilon_i$, $\pm\epsilon_i\pm\epsilon_j$
and $\pm2\delta$, so $\varpi_2^{\mathrm{Sp}(3)}=e_1$.
Let $\mathfrak q=\mathfrak l+\mathfrak u$ be the theta-stable parabolic
whose nilradical consists of the root spaces with $(\alpha,e_1)>0$.
Let $\mathfrak l_{\R}$ be its real Levi algebra.  Then
\[
 \mathfrak l=(\mathfrak l_{\R})_\C\simeq\mathfrak{gl}_1\oplus\mathfrak{so}_7,
 \qquad [\mathfrak l_{\R},\mathfrak l_{\R}]\simeq\mathfrak{so}(4,3),
 \qquad 2\rho(\mathfrak u\cap(\mathfrak p_F)_\C)=6e_1.
\]

For the split Albert algebra, the parameters of Frahm--Labriet are
$r=3$, $d=8$.  Their compact one-dimensional torus is generated by a
rotation in two Jordan coordinates; in the compact $\mathrm{Sp}(3)$
it is conjugate to the direction $\epsilon_1+\epsilon_2=e_1$.
Their restricted roots are $\pm\beta,\pm2\beta$, and their parabolic
is the sum of the zero and positive restricted-root spaces
\cite[Proposition~2.8 and Section~4.2]{FrahmLabriet}.
Thus it is $K_F$-conjugate to $\mathfrak q$ above.  Extending $\beta$
trivially on the semisimple Levi identifies it with $e_1/2$:
the positive $e_1$-coordinates of the $F_4$ roots are $1/2$ and $1$,
with multiplicities eight and seven, respectively.
Writing their integer parameter as $j$ to distinguish it from $k$,
their discrete summand is
\[
 A_{\mathfrak q}(j\beta)\boxtimes\tau^{\mathrm{ds}}_{j+12},
 \qquad j+12\in2\mathbb Z_{>0}.
\]
Here $\tau^{\mathrm{ds}}_{2k}$ has Harish--Chandra module $\sigma_{2k}$.
Consequently
\begin{equation}\label{eq:FL-parameters}
 j=2k-12,\qquad j\beta=(k-6)e_1,\qquad k\geq1.
\end{equation}
The bottom-type shift agrees as well:
$(k-6)e_1+2\rho(\mathfrak u\cap(\mathfrak p_F)_\C)=ke_1$.

\begin{theorem}[Discrete-series lift]\label{thm:lowest-weight-lift}
For every integer $k\geq1$,
\begin{equation}\label{eq:lowest-lift}
 \thetaS_F(\sigma_{2k})\simeq A_{\mathfrak q}((k-6)e_1).
\end{equation}
Its distinguished $K_F$-type is
$V_{\mathrm{Sp}(3)}(k\varpi_2)\boxt\mathbf1$, with multiplicity one.
\end{theorem}

\begin{proof}
By \cite[Theorems~A and~B(2)]{FrahmLabriet} and
\cref{eq:FL-parameters}, the unitary minimal representation contains
a discrete summand whose Harish--Chandra module is
\[
             \sigma_{2k}\boxt A_{\mathfrak q}((k-6)e_1),
\]
and both factors are irreducible.  Orthogonal projection onto this
summand intertwines the subgroup action.  It therefore preserves
smooth vectors and $O(2)\times K_F$-finiteness, and restricts to a
homomorphism
\[
 \PiMin\longrightarrow
             \sigma_{2k}\boxt A_{\mathfrak q}((k-6)e_1).
\]
This homomorphism is nonzero: otherwise the continuous projection
would vanish on the dense ambient Harish--Chandra module and hence
on the entire unitary representation.  Its target is irreducible,
so it is surjective.  By theta adjunction,
$A_{\mathfrak q}((k-6)e_1)$ is an irreducible quotient of
$\ThetaB_F(\sigma_{2k})$.  Part~(3) of \cref{thm:main}, proved in
\cref{sec:pingpong}, gives uniqueness of this quotient and the
asserted multiplicity of its distinguished type.
\end{proof}

\section{Restriction from split \texorpdfstring{$E_6$ to $F_4$}{E6 to F4}}
\label{sec:E6-F4-restriction}

We prove the restriction result used in
\cref{sec:exceptional-value}.  Throughout this section $E_6$ means the
split group $E_{6(6)}$, with compact covering group
$K_6=\mathrm{Sp}(4)$, and $\Pi_6$ is its spherical minimal
Harish--Chandra module.  We compute its quotient by the noncompact
directions of $F_{4(4)}$ on the associated graded level, and then use
unitarity to prove irreducibility of the two summands.

\subsection{The Lagrangian Grassmannian section}

Equip $U=V\oplus W$, with $\dim V=6$ and $\dim W=2$, with the direct
sum symplectic form.  Then
\[
 (\mathfrak p_6)_\C=\Lambda^4_0U,
 \qquad K_F=\mathrm{Sp}(V)\times\mathrm{Sp}(W),
\]
where the subscript $0$ denotes the kernel of symplectic contraction.
The usual ladder filtration of $\Pi_6$ has
\begin{equation}\label{eq:E6-LG-ring}
 R_6:=\gr\Pi_6
   =\bigoplus_{n\geq0}V_{\mathrm{Sp}(4)}(n\varpi_4)
   \simeq\C[\widehat{\operatorname{LG}}(4,U)].
\end{equation}
Here the right side is the coordinate ring of the affine Pl\"ucker cone
over the Lagrangian Grassmannian; its multiplication is the Cartan
product, as in the minimal-representation model of
\cite{BrylinskiKostant}.  In degree one,
\begin{equation}\label{eq:E6-P-T}
 (R_6)_1=P\oplus T,\qquad
 P=\Lambda^3_0V\otimes W=(\mathfrak p_F)_\C,\qquad
 T=\Lambda^2_0V\boxt\mathbf1.
\end{equation}
Thus $\dim P=28$ and $\dim T=14$.  The space $T$ here is the type
$T_1$ of \cref{eq:Pi6-minimal-types}; it is different from the
12-dimensional space $Q$ in Section~\ref{sec:forward-compact-lifts}.

\begin{proposition}[The square-zero section]\label[proposition]{prop:E6-square-zero}
As a graded $K_F$-algebra,
\begin{equation}\label{eq:E6-square-zero}
                  R_6/(P)\simeq\C\oplus T,\qquad T^2=0,
\end{equation}
with $T$ in degree one.  This is a scheme of length $15$ supported at
the origin, with Hilbert series $1+14t$.
\end{proposition}

\begin{proof}
Use the decomposition $\Lambda^4U=A\oplus B\oplus C$ of
\cref{eq:ABC-decomposition}.  Write $\epsilon=x\wedge y$ for the
symplectic bivector of $W$, normalized to have contraction one, and
write an element of $\Lambda^4U$ as
\[
 a+b+c\wedge\epsilon,
 \qquad a\in\Lambda^4V,\quad b\in\Lambda^3V\otimes W,
 \quad c\in\Lambda^2V.
\]
If $\kappa_V$ denotes symplectic contraction on $\Lambda^\bullet V$,
the linear equations defining $\Lambda^4_0U$ are
\begin{equation}\label{eq:E6-isotropy-linear}
 \kappa_Va+c=0,\qquad
 (\kappa_V\otimes1)b=0,\qquad
 \kappa_Vc=0.
\end{equation}
In particular, the middle part is exactly $P$.  Killing its coordinate
functions therefore kills all the $B$-coordinates.

Put $\Omega_V=\sum_{i=1}^3e_i\wedge f_i$.
The contraction map $\kappa_V:\Lambda^4V\to\Lambda^2V$ is an
isomorphism, and
$\kappa_V(\Omega_V\wedge c)=c$ for $c\in\Lambda^2_0V$.
Consequently, after the middle coordinates are killed,
\eqref{eq:E6-isotropy-linear} identifies the remaining linear space with
\[
 \{\,-\Omega_V\wedge c+c\wedge\epsilon:
                      c\in\Lambda^2_0V\,\}.
\]
Restriction of the $A$-coordinates and restriction of the $C$-coordinates
each span the full 14-dimensional space of linear coordinate functions
on this graph.  We identify it with $T$ by symplectic duality.

The ordinary Pl\"ucker relations hold on the Lagrangian Grassmannian.
Since all $B$-coordinates have been killed, the explicit specialization
in \cref{prop:B-linear-section} gives $AC=0$ in $R_6/(P)$.
Both factors span $T$, so every product in $T^2$ is zero.
The ring $R_6$ is generated in degree one, and its quotient has degree-one
space exactly $T$ by \cref{eq:E6-P-T}.  Hence there are no further
linear relations and no surviving terms of degree at least two.
This proves \eqref{eq:E6-square-zero}.
\end{proof}

There is also an exterior-algebra interpretation of this finite scheme.
The wedge-product map identifies
\[
 \frac{\operatorname{Sym}(\Lambda^2V)}
      {(\mathbf S_{(2,2)}V)}
 \simeq\bigoplus_{j\geq0}\Lambda^{2j}V,
\]
with exterior degree $2j$ assigned algebra degree $j$.
Here $\mathbf S_{(2,2)}V$ is the complement of the Pfaffian summand
$\Lambda^4V$ in $\operatorname{Sym}^2(\Lambda^2V)$: killing it makes
products fully alternating in their indices.  Removing the invariant
line means further quotienting by $\Omega_V$.
Since $\Omega_V\wedge-:\Lambda^2V\to\Lambda^4V$ is an isomorphism
in dimension six, this quotient is precisely
$\C\oplus(\Lambda^2V/\C\Omega_V)$ with square-zero positive-degree
part, in agreement with \eqref{eq:E6-square-zero}.

\subsection{The two irreducible summands}

Let $\iota$ be the diagram involution of $E_{6(6)}$ with fixed group
$F_{4(4)}$, preserving the Cartan decomposition.  Twisting by $\iota$
gives again the spherical minimal representation, so there is an
intertwining operator.  Normalize this operator on
$\Pi_6$ to fix the spherical vector.  Its square is then the identity,
and it commutes with $F_{4(4)}$.  On
\eqref{eq:E6-P-T} it acts by $+1$ on $P$ and $-1$ on $T$.

\begin{theorem}[Restriction to $F_{4(4)}$]\label{thm:E6-F4-restriction}
The restriction of $\Pi_6$ is the direct sum of two inequivalent
irreducible unitary Harish--Chandra modules,
\[
                  \Pi_6|_{F_{4(4)}}=\Pi_6^+\oplus\Pi_6^-.
\]
They are the two eigenspaces of $\iota$.  The first is generated by its
spherical vector, and the second by its single copy of
$T=V_{\mathrm{Sp}(3)}(\varpi_2)\boxt\mathbf1$.
\end{theorem}

\begin{proof}
The filtration in \eqref{eq:E6-LG-ring} is $\iota$-stable, and
multiplication by $P$ preserves its two eigenspaces.  By
\cref{prop:E6-square-zero}, the quotient of the positive eigenspace of
$R_6$ by $P$ is the constants, while the corresponding quotient of
the negative eigenspace is $T$.  Graded generation and filtered
Nakayama therefore show that $\Pi_6^+$ is generated over
$U(\mathfrak f_4)$ by the spherical vector and $\Pi_6^-$ by the copy of
$T$ in filtration degree one.  Both eigenspaces are nonzero.

We record the compact branching needed to turn generation into
irreducibility:
\begin{equation}\label{eq:E6-compact-branching}
 V_{\mathrm{Sp}(4)}(n\varpi_4)|_{K_F}
 =\bigoplus_{r=0}^{n}
 V_{\mathrm{Sp}(3)}\bigl((n-r)\varpi_2+r\varpi_3\bigr)
        \boxt\operatorname{Sym}^rW.
\end{equation}
Indeed, symplectic interlacing for the partition $(n,n,n,n)$ leaves
exactly the $\mathrm{Sp}(3)$ partitions $(n,n,r)$, with commuting
$\mathrm{Sp}(1)$-module $\operatorname{Sym}^rW$.
The displayed types are distinct as $(n,r)$ varies.  Thus the whole
restriction is $K_F$-admissible and multiplicity free; in particular,
the trivial type occurs only at $(n,r)=(0,0)$, and $T$ only at $(1,0)$.

The positive-definite invariant Hermitian form on $\Pi_6$ restricts to
both eigenspaces.  Together with $K_F$-admissibility and the finite
generation just proved, this makes each eigenspace a finite direct sum
of irreducible Harish--Chandra modules.  
In such a direct sum, a generating irreducible $K_F$-type of multiplicity
one belongs to only one summand; generation excludes all other summands.
 Both eigenspaces are therefore irreducible.  They are inequivalent because only the
positive eigenspace contains the trivial $K_F$-type.
\end{proof}

\end{document}